\documentclass[11pt]{article}
\usepackage{mathtools}

\usepackage{amsmath,amsfonts,enumerate}
\usepackage{amssymb,amsthm}
\usepackage{calrsfs}
\usepackage{placeins}
\usepackage[margin=1in]{geometry}
\usepackage[T1]{fontenc} \usepackage[utf8]{inputenc} \usepackage{lmodern} \usepackage{graphicx}
\graphicspath{{focm_experiments/}{./}}
\usepackage{biblatex}
\usepackage{csquotes}

\appto{\bibsetup}{\emergencystretch=3em}
\newcommand{\R}{\mathbb{R}}
\newcommand{\inner}[2]{\langle #1, #2 \rangle}
\newcommand{\dG}[2]{d_{G}([#1], [#2])}

\newtheorem{theor}{Theorem}[section]
\newtheorem{lem}[theor]{Lemma}
\newtheorem{prop}[theor]{Proposition}
\newtheorem{cor}[theor]{Corollary}
\newtheorem{defin}[theor]{Definition}
\theoremstyle{definition}

\newtheorem{rems}[theor]{Remarks}
\newtheorem{ex}[theor]{Example}

\newtheorem*{rem*}{Remark}
\DeclareMathOperator{\EX}{\mathbb{E}}

\DeclareMathOperator{\diag}{diag}
\DeclareMathOperator{\Cov}{Cov}

\title{Distortion of Random Max Filters over Reflection Groups}
\author{{\bf Efstratios Tsoukanis} \\
	Institute of Mathematical Sciences \\ Claremont Graduate University \\
	Claremont, CA, USA \\ email: efstratios.tsoukanis@cgu.edu}

\date{\today}
\begin{document}

	\maketitle

\begin{abstract}
	A max filter with Gaussian templates $a_1,\dots,a_D$ is a 
	$G$-invariant map $\R^d\to\R^D$, where $G$ is a finite reflection group. As
	$D\to\infty$ its distortion converges a.s to the condition number of
	$M_G^{1/2}$, where $M_G$ is the second-moment matrix of a standard Gaussian
	vector moved into a fundamental chamber of $G$. We compute this limit exactly
	for the sign and dihedral groups, and show it grows as $\Theta(\sqrt d)$,
	$\Theta(d^{3/2})$ and $\Theta(k)$ for the sign, symmetric  and dihedral groups of order $2k$. The
	symmetric case bound uses a lower bound on the covariance of order
	statistics. For any finite $G$ we also show that the infimum defining the optimal lower Lipschitz
	constant of a piecewise-linear $G$-invariant map is a minimum.
\end{abstract}

	\section{Introduction}

	In many problems in machine learning and signal processing, the data is known
	 up to some group action. Some examples are molecular
	conformation spaces~\cite{dym2022low}, crystallographic
	reconstruction~\cite{edidin2023generic}, and permutation-invariant deep
	networks~\cite{deepsets,balan2025permutation}. In these problems we need an invariant map $f \colon \R^d \to \R^m$, $f(g \cdot x) = f(x)$, that is injective on the
	orbit space $\R^d/G$, so that distinct orbits are separated, and stable. For a large class of maps injectivity implies stability~\cite{freeman2024stability,sverdlov2024fsw}. For example, if $f$ is continuous, positively homogeneous and piecewise linear, then injectivity of $f$ implies that $f$ is bi-Lipschitz (Theorem~\ref{thm:main} below). This class contains max filters, sorting-based and
	coorbit embeddings, and bias free invariant ReLU networks.

	\subsection*{Notation}
	Let $G$ be a finite group acting orthogonally on $V \cong \R^d$, with $g$ acting
	via the orthogonal map $U_g$. For $x\in V$ we write $[x]:=\{U_gx:g\in G\}$ for
	its orbit. The natural distance on the quotient $V/G$ is the orbit distance:
	\[
	\dG{x}{y} \;:=\; \min_{g \in G} \|x - U_g y\|_2.
	\]
	With a slight abuse of notation we write $d_G(x,y)$ for $\dG{x}{y}$.
	Given a $G$-invariant map $f \colon V \to \R^m$, we say $f$ is bi-Lipschitz
		 if there exist constants $0 < a \leq b < \infty$ such that
	\begin{equation}\label{eq:biLip}
		a\,d_G(x,y) \;\leq\; \|f(x) - f(y)\|_2 \;\leq\; b\,d_G(x,y)
		\qquad \forall\, x, y \in V.
	\end{equation}
	The optimal lower and upper Lipschitz constants of $f$ are
	\[
	a(f):=\inf_{[x]\neq[y]}\frac{\|f(x)-f(y)\|_2}{d_G(x,y)},\qquad
	b(f):=\sup_{[x]\neq[y]}\frac{\|f(x)-f(y)\|_2}{d_G(x,y)},
	\]
	with values in $[0,\infty)$.  In this case $a(f)$ is the largest $a$ and $b(f)$ is the smallest $b$
	in~\eqref{eq:biLip}, and we call the ratio $\kappa(f):=b(f)/a(f)$ the
	\emph{distortion} of $f$. Our goal is to compute or
	estimate $a(f)$, $b(f)$ and $\kappa(f)$. For a linear map $A$ we write
	$\sigma_{\max}(A):=\max_{\|z\|_2=1}\|Az\|_2$ and
	$\sigma_{\min}(A):=\min_{\|z\|_2=1}\|Az\|_2$.

	\subsection*{Class of Maps}

	\begin{defin}\label{piecelinear}
		A function $f \colon \R^d \to \R^m$ is called piecewise-linear $G$-invariant
			positively homogeneous if:
		\begin{enumerate}[(i)]
			\item There exist finitely many affine maps $f_1,\dots,f_N$ from $\R^d$ to
			$\R^m$ and a map $s:\R^d\to\{1,\dots,N\}$ such that
			$f(x)=f_{s(x)}(x)$ for every $x \in \R^d$;
			\item $f$ is continuous on $\R^d$;
			\item $f(U_g x) = f(x)$ for all $g \in G$, $x \in \R^d$;
			\item $f(\lambda x) = \lambda f(x)$ for all $\lambda \geq 0$, $x \in \R^d$.
		\end{enumerate}
		 We call each $f_i$ a component of $f$.
	\end{defin}

	\begin{rem*}
		Under (iv), the components can be taken linear. Write $f_i=f_i^{0}+b_i$ with
		$f_i^{0}$ linear and $b_i\in\R^m$, and fix $x\in\R^d$. Some index $i$ satisfies $s(nx)=i$
		for infinitely many $n\in\mathbb N$. For these $n$, (iv) gives
		$nf(x)=f(nx)=nf_i^{0}(x)+b_i$, so $b_i=n\bigl(f(x)-f_i^{0}(x)\bigr)$. Since $b_i$ does not
		depend on $n$, we get $f(x)=f_i^{0}(x)$. So, $f(x)\in\{f_1^{0}(x),\dots,f_N^{0}(x)\}$
		for every $x$, and we can replace each $f_i$ by $f_i^{0}$, with a new map that
		sends $x$ to an index $i$ with $f(x)=f_i^{0}(x)$.
	\end{rem*}

	These properties give $f$ a polyhedral structure. Finitely many hyperplanes
	decompose $\R^d$ into open convex polyhedral cones, and on each open cone $f$ agrees
	with a single linear component (Proposition~\ref{prop:conic-decomp}).

	\subsection*{Injectivity Implies Stability}

	\begin{theor}\label{thm:main}
		Let $G$ be a finite group acting orthogonally on $V \cong \R^d$, and let
		$f \colon V \to \R^m$ be as in Definition~\ref{piecelinear}. If $f$ is injective
		on the quotient space $V/G$, then there exist constants $0 < a \leq b < \infty$
		such that
		\[
		a\,d_G(x,y) \;\leq\; \|f(x) - f(y)\|_2 \;\leq\; b\,d_G(x,y)
		\]
		for all $x, y \in V$.
	\end{theor}

	Theorem~\ref{thm:main} follows from the ~\cite[Lem.~3.4]{sverdlov2024fsw}, which we state and apply below. Earlier special cases include multisets of \(n\) points in \(\R^{d'}\) under row permutations, by the proof of~\cite[Thm.~4.2]{amir2025fourier}, as well as coorbit and sorting embeddings~\cite{balan2023stab,balan2025permutation,balan2024sorting}. The general result was announced in~\cite{freeman2024stability}. Here, a function on a compact convex polytope \(P\) is called piecewise linear if \(P\) is a finite union of compact convex polytopes on each of which the function is affine.

\begin{lem}[Sverdlov et al., polytope formulation]
	\label{lem:sdda}
	Let $P\subset\R^k$ be a compact convex polytope, and let
	$F,Q\colon P\to[0,\infty)$ be piecewise linear. Suppose that,
	for every $z\in P$, $F(z)=0$ if and only if $Q(z)=0$.
	Then there are constants $c,C>0$ such that
	\[
	c\,Q(z)\le F(z)\le C\,Q(z)
	\qquad\text{for all }z\in P.
	\]
\end{lem}

	This is Lemma~3.4 of~\cite{sverdlov2024fsw}, which is stated there for a compact
	polygon $P$. Note that we only need the case where $P$ is convex. In~\cite{sverdlov2024fsw} it is used to prove that the
	FSW-GNN architecture is bi-Lipschitz. It is also noted there that the lemma implies
	earlier bi-Lipschitz results, among them those
	of~\cite{balan2023stab,balan2024sorting,balan2025permutation}, because the maps are positively homogeneous. These earlier proofs use homogeneity to reduce to pairs of points in a bounded set, and then
	analyze locally.

\begin{proof}[Proof of Theorem~\ref{thm:main} as a corollary of Lemma~\ref{lem:sdda}]
	Let $P:=[-1,1]^{2d}\subset V\times V$, and for
	$(x,y)\in V\times V$ put
	\[
	F(x,y):=\|f(x)-f(y)\|_1,\qquad
	Q(x,y):=\min_{g\in G}\|x-U_gy\|_1.
	\]
	We first show that $F$ and $Q$ are piecewise linear on $P$.
	
	By the remark after Definition~\ref{piecelinear} and
	Proposition~\ref{prop:conic-decomp}, the finitely many open
	cones of that proposition have dense union. Their closures
	therefore cover $V$, and by continuity $f$ is linear on
	each closure.
	
	On the product of any two such closed cones, the map
	$(x,y)\mapsto f(x)-f(y)$ is linear.
	Intersecting these products with $P$ and dividing according
	to the signs of the coordinates of $f(x)-f(y)$, we obtain
	finitely many polytopes on which $F$ is affine.
	
	For $Q$, we divide $P$ according to the signs of the
	coordinates of $x-U_gy$, for every $g\in G$.
	On each resulting polytope, every function
	$\|x-U_gy\|_1$ is affine.
	We divide once more according to which $g$ gives the
	minimum. Then $Q$ is affine on each final piece.
	All these divisions are taken along hyperplanes, so the pieces
	remain compact convex polytopes. Taking a common
	subdivision, both $F$ and $Q$ are piecewise linear on $P$.
	
	Next, $Q(x,y)=0$ exactly when $x=U_gy$ for some
	$g\in G$, while $F(x,y)=0$ exactly when $f(x)=f(y)$.
	By $G$-invariance and injectivity on $V/G$, both conditions
	are equivalent to $[x]=[y]$.
	Thus Lemma~\ref{lem:sdda} gives constants $c,C>0$ with
	$c\,Q\le F\le C\,Q$ on $P$.
	
	By (iv) of Definition~\ref{piecelinear} and the linearity
	of the action,
	\[
	F(tx,ty)=tF(x,y),\qquad
	Q(tx,ty)=tQ(x,y),\qquad t\ge0.
	\]
	Since every pair $(x,y)$ is a nonnegative multiple of
	a point of $P$, the bounds $c\,Q\le F\le C\,Q$
	hold on all of $V\times V$.
	
	Finally, for $v\in\R^k$ we have
	$\|v\|_2\le\|v\|_1\le\sqrt{k}\,\|v\|_2$.
	Taking the minimum over $g\in G$ gives
	$d_G(x,y)\le Q(x,y)\le\sqrt d\,d_G(x,y)$.
	Hence
	\[
	\|f(x)-f(y)\|_2
	\ge \frac{F(x,y)}{\sqrt m}
	\ge \frac{c}{\sqrt m}\,Q(x,y)
	\ge \frac{c}{\sqrt m}\,d_G(x,y)
	\]
	and
	\[
	\|f(x)-f(y)\|_2
	\le F(x,y)
	\le C\,Q(x,y)
	\le C\sqrt d\,d_G(x,y).
	\]
	Thus the theorem holds with $a:=c/\sqrt m$
	and $b:=C\sqrt d$. 
\end{proof}

	Note that this argument gives no information about the size of the constants. In
	Section~\ref{sec:formula} we give a different proof of Theorem~\ref{thm:main}, which
	goes through the optimal lower constant (Corollary~\ref{cor:main-from-formula}). 

	\subsection*{Main Results}

	\noindent\textbf{1. Random max filters over reflection groups.} For templates
	$A=(a_1,\dots,a_D)$, $a_j\in\R^d$, the \emph{max filter} \cite{cahill2025group} is the map $f_A\colon\R^d\to\R^D$,
	$f_A(x)_j:=\max_{g\in G}\inner{a_j}{U_gx}$, $j\in[D]$. Now let $G$ be generated by reflections, and let $\Lambda_G$ be the chamber chart of Section~\ref{subsec:chamber}. In~\cite[Lem.~9]{mixonpacker} it is shown that $f_A$ has the
	same optimal constants as the $D\times d$ matrix whose rows are the
	templates moved into the chamber by $\Lambda_G$. If $a_1,a_2,\dots$ are i.i.d.\
	$N(0,I_d)$, then for fixed $d$, almost surely,
	\[
	\lim_{D\to\infty}\kappa(f_A)=\kappa_G=\Bigl(\frac{\lambda_{\max}(M_G)}{\lambda_{\min}(M_G)}\Bigr)^{1/2},
	\qquad
	M_G=\EX\bigl[\Lambda_G(\xi)\Lambda_G(\xi)^{\!\top}\bigr],\quad \xi\sim N(0,I_d)
	\]
	(Definition~\ref{def:MG} and Theorem~\ref{thm:spectral-limit}). We  estimate $\kappa_G$ for the following groups
	(Theorems~\ref{thm:exact} and~\ref{thm:growth}):
	\begin{center}
		\renewcommand{\arraystretch}{1.3}
		\begin{tabular}{ll}
			\hline
			Group & Limiting distortion $\kappa_G$ \\
			\hline
			sign group $\{\pm1\}^d$, $d\ge2$ & $\sqrt{1+2d/(\pi-2)}=(1+o(1))\sqrt{2/(\pi-2)}\,\sqrt d$ \\
			$S_d$ on $\R^d$ or on $\mathbf 1^{\perp}$ & $\Theta(d^{3/2})$, at most $\sqrt{\tfrac2\pi\,d(d+1)(d+2)}$ \\
			$B_d$ and $D_d$ & $\Theta(d^{3/2})$, at most $\sqrt{\tfrac8\pi\,d(d+1)(d+2)}$ \\
			dihedral $I_2(k)$ & $\bigl(\tfrac{1+\sin w/w}{1-\sin w/w}\bigr)^{1/2}=(1+O(k^{-2}))\tfrac{2\sqrt3}{\pi}k$, $w=\tfrac\pi k$ \\
			\hline
		\end{tabular}
	\end{center}

	\smallskip
	\noindent\textbf{2. An estimate for order statistics.}
	We prove a lower bound of order $d^{-2}$ for the smallest
	eigenvalue of the covariance matrix of order statistics
	(Lemma~\ref{lem:orderstat-floor}), and a matching upper bound
	under an additional assumption on the density
	(Lemma~\ref{lem:central-spacing}).
	The rate and constant of the lower bound agree with a classical
	heuristic approximation~\cite{stephens1990asymptotic},
	\cite[Ch.~4]{david2003order}.
	For the uniform distribution, the bound is asymptotically
	sharp~\cite{durbin1972components}.
	
	\smallskip
	\noindent\textbf{3. Optimal constants of piecewise-linear invariants.}
	For any finite group $G$, we express $a(f)$ as the minimum of
	finitely many local problems and prove that each infimum is
	a minimum (Proposition~\ref{prop:mu-attained}).
	Consequently, $a(f)$ is achieved at a pair of distinct orbits
	(Corollary~\ref{cor:main-from-formula}).
	For max filters, this was known in special cases
	\cite[Lem.~27]{mixonqaddura},
	\cite[Lem.~9]{mixonpacker}.
	We also express $b(f)$ as the largest operator norm of the
	linear maps that agree with $f$ on an open cone
	(Proposition~\ref{thm:formula}(iii)).
	This was sketched in~\cite{ghoreishi2024maxfilter}, as reported
	in~\cite[Sec.~3.1]{mixonqaddura}, and proved for max filters
	in~\cite[Thm.~14]{mixonqaddura}.

\subsection*{Relation to Prior Work}

 Whether the quotient of a Hilbert space by a group of its automorphisms admits a
bi-Lipschitz embedding has been studied
in~\cite{cahill2024towards}. Separating invariants for finite and reflection groups
are studied in~\cite{DUFRESNE20091979,compuinvar,dadok1985polar,edidin2025orbit}, and
bi-Lipschitz embeddings of orbit spaces are constructed
in~\cite{eriksson2018quantitative,weighill2023coarse,derksen2024bilipschitz,blumsmith2025distortion}.
For finite groups, stability has been shown for max
filtering~\cite{cahill2025group,mixonpacker,mixonqaddura,cahill2026approx},
and sorting-based embeddings~\cite{balan2023stab,balan2023inj,balan2024sorting}, by
semialgebraic methods~\cite{bendory2025phase,bendory2025generalized,amir2025stability},
and through tropical polynomials~\cite{derksen2025tropical,kubo2026stable}.
For compact groups, local bi-Lipschitz stability of max filters on
orbits is studied in~\cite{qaddura2025max}, with applications
to weighted phase retrieval and cryo-EM.
On the other hand, in~\cite{davidson2024holder} it is shown that sum-pooling with ReLU
activation is, in expectation over the parameters and on bounded sets, lower H\"older but not
lower Lipschitz.

Less is known about the optimal constants. For random templates, the known bounds on the
distortion of max filters grow with
$|G|$~\cite[Thm.~32]{mixonqaddura},~\cite[Thm.~18]{cahill2025group}, and they are not believed
to be tight~\cite{cahill2026approx}. Note that the bound of~\cite[Thm.~18]{cahill2025group} is
stated for templates drawn uniformly from the sphere, and it also depends on $d$. For
phase retrieval the exact limiting distortion is
determined in~\cite{xia2024stability}. For sorting-based embeddings of multisets of $n$
points in $\R^{d'}$ with $d'\ge2$, it is shown in~\cite{dym2026quant} that the distortion is
at least of order $\sqrt n$ and, with high
probability, at most of order $n^2$ for $D\gtrsim n^4d'$ Gaussian directions.

For reflection groups the orbit space embeds isometrically, so its Euclidean distortion
satisfies $c_2(\R^d/G)=1$~\cite{blumsmith2025distortion,cahill2026approx}. Moreover, a
max-filter bank with $d$ templates, followed by a linear map, achieves this distortion
exactly~\cite[Thm.~3]{cahill2026approx}.

\subsection*{Organization}

\begin{itemize}
	\item In Section~\ref{sec:random-maxfilter} we present our results on random max filters
	over reflection groups.
	\item In Section~\ref{app:orderstat} we prove two lemmas on order statistics.
	\item In Section~\ref{sec:proofs} we prove
	Theorem~\ref{thm:growth}, we treat the sign group when $d\to\infty$ and
	$D/d\to\infty$, and we give a version of Theorem~\ref{thm:spectral-limit} for a finite
	number of templates.
	\item In Section~\ref{sec:formula} we consider any finite group. We describe the structure
	of the maps in our class, and we write the optimal lower constant as a minimum of finitely
	many local problems. Then we show that the infimum in each of them is a minimum, and we give
	a second proof of Theorem~\ref{thm:main}. 
	\item In Section~\ref{sec:experiments} we test the growth laws numerically.
\end{itemize}

\section{Random Max Filters over Reflection Groups}
\label{sec:random-maxfilter}

\subsection{The Chamber Chart}\label{subsec:chamber}

Let $G\le O(d)$ be a finite group generated by orthogonal reflections. Fix a
system of simple roots and let
\[
C:=\bigl\{x\in\R^d:\ \inner{x}{\alpha}\ge0\ \text{for every simple root }\alpha\bigr\}
\]
be the chosen fundamental chamber. Such a system
exists~\cite[Thm.~1.3]{humphreys}. Every orbit has  exactly one
point in $C$  ~\cite[Thm.~1.12]{humphreys}. We call the map $\Lambda_G:\R^d\to C$ that
sends $x$ to the  point of $[x]\cap C$ the chamber chart. 
$|x|:=(|x_1|,\dots,|x_d|)$, and ${\downarrow}x$ is the vector of the coordinates of
$x$ in decreasing order. For the sign group $G=\{\pm1\}^d$ let
$C=\R_{\ge0}^d$, and then $\Lambda_G(x)=|x|$. For $G=S_d$ let
$C=\{x_1\ge\cdots\ge x_d\}$, and then $\Lambda_G(x)={\downarrow}x$. For the group
$B_d$ of signed permutations one may take $C=\{x_1\ge\cdots\ge x_d\ge0\}$, and then
$\Lambda_G(x)={\downarrow}|x|$. For the group $D_d$ of signed permutations with an even
number of sign changes let $C=\{x_1\ge\cdots\ge x_{d-1}\ge|x_d|\}$, and then
$\Lambda_G(x)=(y_1,\dots,y_{d-1},\varsigma(x)\,y_d)$, where $y:={\downarrow}|x|$ and
$\varsigma(x):=\prod_i\operatorname{sign}(x_i)$. For $S_d$ acting on $\mathbf 1^{\perp}$ let $C=\{x_1\ge\cdots\ge x_d\}\cap\mathbf 1^{\perp}$, and  $\Lambda_G(x)={\downarrow}x$.
For $k\ge2$, the dihedral group $I_2(k)$ of order $2k$ acts on $\R^2$. It is generated by
the reflections in the lines spanned by $(1,0)$ and by $(\cos\frac\pi k,\sin\frac\pi k)$.
Let $C:=\{r(\cos\theta,\sin\theta): r\ge0,\ 0\le\theta\le\pi/k\}$, the
sector between these lines.

In the following lemma we collect the properties of $\Lambda_G$ that we use. 
Parts (b) and (c) are known, and we do not prove them here. Part (c) is \cite[Lem.~8]{mixonpacker} and part (b)
is \cite[eq.~(2)]{mixonpacker}. 
\begin{lem}
	\label{lem:chamber-chart}
	Let $G\le O(d)$ be a finite reflection group with closed chamber $C$ and
	chart $\Lambda_G$. Then:
	\begin{enumerate}[(a)]
		\item $\Lambda_G$ is positively homogeneous and
		$\|\Lambda_G(x)\|_2=\|x\|_2$ for all $x$;
		\item $\displaystyle\max_{g\in G}\,\inner{a}{U_gx}
		=\bigl\langle\Lambda_G(a),\Lambda_G(x)\bigr\rangle$
		for all $a,x\in\R^d$;
		\item $d_G(x,y)=\bigl\|\Lambda_G(x)-\Lambda_G(y)\bigr\|_2$
		for all $x,y\in\R^d$;
		\item $C-C=\R^d$.
	\end{enumerate}
\end{lem}

\begin{proof}[Proof of (a) and (d)]
	(a) The action commutes with scaling and $C$ is a cone, so
	$\Lambda_G(\lambda x)=\lambda\Lambda_G(x)$ for $\lambda\ge0$. Moreover, each $U_g$ is
	orthogonal, and therefore $\|\Lambda_G(x)\|_2=\|x\|_2$.

	(d) Let $v\in\R^d$. Let $x_0$ an interior point on $C$,	so $x_0+tv\in C$ for some $t>0$. Then $v=t^{-1}(x_0+tv)-t^{-1}x_0$ is a difference
	of two elements of the cone $C$.
\end{proof}

\subsection{Exact Constants for Arbitrary Templates}

The following theorem is~\cite[Lem.~9]{mixonpacker}, restated in
our notation, with the rank criterion added.

\begin{theor}[{\cite[Lem.~9]{mixonpacker}}]
	\label{thm:reflection-constants}
	Let $G\le O(d)$ be a finite reflection group, let $a_1,\dots,a_D\in\R^d$,
	and let $R\in\R^{D\times d}$ be the matrix with rows
	$\Lambda_G(a_j)^{\!\top}$. Then $f_A=R\,\Lambda_G$,
	\[
	a(f_A)=\sigma_{\min}(R)\qquad\text{and}\qquad
	b(f_A)=\sigma_{\max}(R).
	\]
	Moreover, the following are equivalent: $f_A$ separates orbits, $f_A$ is
	bi-Lipschitz, and $\operatorname{rank}R=d$. In that case
	$\kappa(f_A)=\kappa(R)$, the condition number of $R$.
\end{theor}

\begin{proof}
	The $f_A=R\,\Lambda_G$ follows from Lemma~\ref{lem:chamber-chart}(b). By
	Lemma~\ref{lem:chamber-chart}(c), for $[x]\neq[y]$,
	\[
	\frac{\|f_A(x)-f_A(y)\|_2}{d_G(x,y)}
	=\frac{\bigl\|R\bigl(\Lambda_G(x)-\Lambda_G(y)\bigr)\bigr\|_2}
	{\bigl\|\Lambda_G(x)-\Lambda_G(y)\bigr\|_2},
	\]
	and as $(x,y)$ ranges over all pairs with $[x]\neq[y]$, the difference
	$\Lambda_G(x)-\Lambda_G(y)$ ranges over all of
	$(C-C)\setminus\{0\}=\R^d\setminus\{0\}$ by
	Lemma~\ref{lem:chamber-chart}(d). Hence $a(f_A)=\sigma_{\min}(R)$ and
	$b(f_A)=\sigma_{\max}(R)$.

	If $\operatorname{rank}R=d$, then $a(f_A)>0$ and $f_A$ is bi-Lipschitz so injective. On the other hand, if $\operatorname{rank}R<d$, we pick $v\neq0$ with $Rv=0$ and we write
	$v=z-z'$ with $z,z'\in C$. Then $z\neq z'$ are representatives of distinct
	orbits with $f_A$ values $Rz=Rz'$, so $f_A$ is not injective.
\end{proof}

\subsection{Random Templates: A Spectral Reduction}

\begin{defin}\label{def:MG}
	For a finite reflection group $G\le O(d)$ let $\xi\sim N(0,I_d)$, and set
	\[
	\mu_G:=\EX[\Lambda_G(\xi)],\qquad
	\Sigma_G:=\Cov\bigl(\Lambda_G(\xi)\bigr),\qquad
	M_G:=\EX\bigl[\Lambda_G(\xi)\,\Lambda_G(\xi)^{\!\top}\bigr]=\Sigma_G+\mu_G^{}\mu_G^{\!\top}.
	\]
	By Lemma~\ref{lem:chamber-chart}(a),
	$\operatorname{tr}M_G=\EX\|\Lambda_G(\xi)\|_2^2=\EX\|\xi\|_2^2=d$.
\end{defin}

The following theorem is a consequence of Theorem~\ref{thm:reflection-constants}.

\begin{theor}
	\label{thm:spectral-limit}
	Let $G\le O(d)$ be a finite reflection group. Let
	$a_1,a_2,\dots$ be i.i.d.\ $N(0,I_d)$, and for each $D$ let $f_A$ be the max filter
	with templates $a_1,\dots,a_D$. Then $\Sigma_G\succ0$ and $M_G\succ0$. Moreover, almost surely $a(f_A)>0$ for every $D\ge d$, and
	\[
	\kappa(f_A)\ \longrightarrow\
	\kappa_G:=\Bigl(\frac{\lambda_{\max}(M_G)}{\lambda_{\min}(M_G)}\Bigr)^{1/2}
	\qquad\text{as }D\to\infty.
	\]
\end{theor}

\begin{proof}
	We write $Z:=\Lambda_G(\xi)$, with $\xi$ as in Definition~\ref{def:MG}. For each $D$, let
	$R$ be the matrix of Theorem~\ref{thm:reflection-constants} for the templates
	$a_1,\dots,a_D$. Its $j$-th row is $Z_j^{\!\top}$, where $Z_j:=\Lambda_G(a_j)$.

	Let $B\subseteq\R^d$ be a Borel set. If $Z\in B$, then $\xi\in U_gB$ for some $g\in G$. Since the law of
	$\xi$ is invariant under $O(d)$,
	\[
	\mathbb{P}(Z\in B)\ \le\ \sum_{g\in G}\mathbb{P}(\xi\in U_gB)=|G|\,\mathbb{P}(\xi\in B).
	\]

Now let $v\neq0$, and let $H:=\{z\in\R^d:\inner{v}{z}=\inner{v}{\mu_G}\}$. The bound above with
$B=H$ gives $\mathbb{P}(Z\in H)=0$. So $\inner{v}{Z}\neq\inner{v}{\mu_G}$ almost surely, and
hence $v^{\!\top}\Sigma_Gv=\operatorname{Var}\inner{v}{Z}>0$. Therefore $\Sigma_G\succ0$, and
$M_G=\Sigma_G+\mu_G\mu_G^{\!\top}\succeq\Sigma_G\succ0$.

	Next, by the bound above, $Z$ lies in a fixed hyperplane with probability $0$.  So, the first $d$ rows of $R$ are linearly
	independent almost surely. These rows belong to $R$ for every $D\ge d$, so almost surely
	$\operatorname{rank}R=d$ for every $D\ge d$, and $a(f_A)=\sigma_{\min}(R)>0$ by
	Theorem~\ref{thm:reflection-constants}.

	Finally, $\EX\|Z\|_2^2=\EX\|\xi\|_2^2=d$ by Lemma~\ref{lem:chamber-chart}(a), so every
	entry of $ZZ^{\!\top}$ is integrable. By the strong law of large numbers, applied to each
	entry, $R^{\!\top}R/D=\frac1D\sum_{j=1}^DZ_jZ_j^{\!\top}\to M_G$ almost surely. Moreover, eigenvalues depend continuously on
	the matrix, and therefore $\lambda_{\max}(R^{\!\top}R/D)\to\lambda_{\max}(M_G)$ and
	$\lambda_{\min}(R^{\!\top}R/D)\to\lambda_{\min}(M_G)>0$. Since $\operatorname{rank}R=d$,
	by Theorem~\ref{thm:reflection-constants} we have that 
	\[
	\kappa(f_A)=\frac{\sigma_{\max}(R)}{\sigma_{\min}(R)}
	=\Bigl(\frac{\lambda_{\max}(R^{\!\top}R/D)}{\lambda_{\min}(R^{\!\top}R/D)}\Bigr)^{1/2}
	\ \longrightarrow\ \kappa_G 
	\]
	almost surely.
\end{proof}

\begin{theor}
	\label{thm:exact}
	Let $M_G$ defined in ~\ref{def:MG}, with the chambers in ~\ref{subsec:chamber}, and  $\kappa_G$ be as in Theorem~\ref{thm:spectral-limit}.
	\begin{enumerate}[(a)]
		\item \emph{(Sign group $G=\{\pm1\}^d$.)} For every $d$,
		\[
		M_G=\Bigl(1-\tfrac2\pi\Bigr)I_d+\tfrac2\pi\,\mathbf 1\mathbf 1^{\!\top}.
		\]
		Its eigenvalues are $1-\tfrac2\pi+\tfrac{2d}\pi$, with eigenvector $\mathbf 1$, and,
		if $d\ge2$, $1-\tfrac2\pi$, with multiplicity $d-1$. Hence, for $d\ge2$,
		\[
		\kappa_G=\sqrt{1+\frac{2d}{\pi-2}}
		=\bigl(1+o(1)\bigr)\sqrt{\tfrac{2}{\pi-2}}\;\sqrt d .
		\]

		\item \emph{(Dihedral group $G=I_2(k)$ acting on $\R^2$, $k\ge2$.)}
		Let $w:=\pi/k$. For every $k$, the eigenvalues of $M_G$ are $1\pm\frac{\sin w}{w}$, so
		\[
		\kappa_G=\Bigl(\frac{1+\tfrac{\sin w}{w}}{1-\tfrac{\sin w}{w}}\Bigr)^{1/2}
		=\bigl(1+O(k^{-2})\bigr)\,\frac{2\sqrt3}{\pi}\,k .
		\]
	\end{enumerate}
\end{theor}

\begin{proof}
	(a) Here $\Lambda_G(\xi)=|\xi|$, whose coordinates are i.i.d.\ half-normal, with
	$\EX|\xi_1|^2=1$ and $\EX|\xi_1|=\sqrt{2/\pi}$. Therefore the diagonal entries of $M_G$ are equal to
	$1$, and the off-diagonal entries are equal to $(\EX|\xi_1|)^2=2/\pi$. For $d\ge2$ both eigenvalues occur, and so
	$\kappa_G^2=\lambda_{\max}(M_G)/\lambda_{\min}(M_G)=1+\frac{2d/\pi}{1-2/\pi}=1+\frac{2d}{\pi-2}$.

	(b) The lines of $I_2(k)$ divide the plane into $2k$ sectors of angle
	$w$, and each sector is mapped onto $C$ by one group element. A rotation of
	coordinates conjugates $M_G$ by an orthogonal matrix and does not change its
	eigenvalues. Thus we may assume that $C$ is the sector of polar angles $[0,w]$. We write
	$\xi=r(\cos\theta,\sin\theta)$ with $r\ge0$. By rotation invariance of the Gaussian,
	$\theta$ is uniform on $[0,2\pi)$ and independent of $r$, and $\EX r^2=2$. Each group
	element is a rotation or a reflection, so it preserves $r$ and maps the angles of a
	sector onto $[0,w]$ by $\theta\mapsto\pm\theta+c$ modulo $2\pi$. Hence
	$\Lambda_G(\xi)=r(\cos\theta^*,\sin\theta^*)$, where $\theta^*\in[0,w]$ is the image of
	$\theta$. Each sector has probability $1/(2k)$
	so $\theta^*$ is uniform on $[0,w]$.Using $\cos^2\theta=\frac{1+\cos2\theta}2$,
	$\sin^2\theta=\frac{1-\cos2\theta}2$ and $\cos\theta\sin\theta=\frac{\sin2\theta}2$, we
	get
	\[
	M_G=\EX[r^2]\,\EX\begin{pmatrix}\cos^2\theta^*&\cos\theta^*\sin\theta^*\\
	\cos\theta^*\sin\theta^*&\sin^2\theta^*\end{pmatrix}
	=I_2+\begin{pmatrix}p&q\\q&-p\end{pmatrix},
	\qquad p:=\EX\cos2\theta^*,\quad q:=\EX\sin2\theta^*.
	\]
	Since $\theta^*$ is uniform on $[0,w]$,
	\[
	p=\frac1w\int_0^w\cos2\theta\,d\theta=\frac{\sin w\cos w}{w},
	\qquad
	q=\frac1w\int_0^w\sin2\theta\,d\theta=\frac{\sin^2w}{w}.
	\]
	Thus $p^2+q^2=\sin^2w/w^2$. The second matrix has trace $0$ and determinant
	$-(p^2+q^2)$, so its eigenvalues are $\pm\sin w/w$, because $\sin w>0$. Therefore the eigenvalues of $M_G$ are $1\pm\sin w/w$. Finally, $1-\sin w/w=\tfrac{w^2}{6}\bigl(1+O(w^2)\bigr)$ and
	$1+\sin w/w=2\bigl(1+O(w^2)\bigr)$, so
	$\kappa_G=(1+O(w^2))\,\sqrt{12}/w=(1+O(k^{-2}))\,\tfrac{2\sqrt3}{\pi}k$.
\end{proof}

\begin{theor}
	\label{thm:growth}
	Let $M_G$ and $\kappa_G$ be as in Theorem~\ref{thm:exact}.
	\begin{enumerate}[(a)]
		\item \emph{(Symmetric group $G=S_d$ acting on $\R^d$ by permuting the coordinates.)} For an absolute constant $C_1$,
		\[
		\begin{aligned}
			(1-o(1))\,d\ &\le\ \lambda_{\max}(M_G)\ \le\ d,\\
			\frac{\pi}{2(d+1)(d+2)}\ &\le\ \lambda_{\min}(M_G)\ \le\ \frac{C_1}{d^{2}} .
		\end{aligned}
		\]
		Hence $\kappa_G=\Theta(d^{3/2})$, and
		\[
		\kappa_G\ \le\ \sqrt{\tfrac{2}{\pi}\,d(d+1)(d+2)}
		\ =\ \bigl(\sqrt{2/\pi}+o(1)\bigr)\,d^{3/2}.
		\]
		 For $d\ge2$, let $G'$ be $S_d$ acting on $\mathbf 1^{\perp}\cong\R^{d-1}$. Then
		$\lambda_{\max}(M_{G'})=\lambda_{\max}(M_G)$ and $\lambda_{\min}(M_{G'})=\lambda_{\min}(M_G)$.
		In particular $\kappa_{G'}=\kappa_G$, and the bounds above also hold for $G'$. 
		\item \emph{(Signed permutations: $G=B_d$ and $G=D_d$.)} In both
		cases, for an absolute constant $C_1$,
		\[
		(1-o(1))\,d\ \le\ \lambda_{\max}(M_G)\ \le\ d,
		\qquad
		\frac{\pi}{8(d+1)(d+2)}\ \le\ \lambda_{\min}(M_G)\ \le\ \frac{C_1}{d^{2}}.
		\]
		Hence $\kappa_G=\Theta(d^{3/2})$, and
		\[
		\kappa_G\ \le\ \sqrt{\tfrac{8}{\pi}\,d(d+1)(d+2)}
		\ =\ \bigl(\sqrt{8/\pi}+o(1)\bigr)\,d^{3/2}.
		\]
	\end{enumerate}
\end{theor}

We prove Theorem~\ref{thm:growth} in Section~\ref{app:bc}, using the lemmas on order statistics of
Section~\ref{app:orderstat}.

\begin{rems}\leavevmode
	\label{rems:reflection}
	\begin{enumerate}[(i)]
		\item Reflection groups have Voronoi characteristic $1$~\cite[Lem.~41]{mixonqaddura}. Therefore, as $D\to\infty$, the distortion bound of~\cite[Thm.~32]{mixonqaddura} for random templates
		tends to a constant times $|G|\log^{1/2}(e|G|)$. For $S_d$, $B_d$ and $D_d$ this
		is superexponential in $d$, while Theorem~\ref{thm:growth} gives $\Theta(d^{3/2})$. For
		$S_d$, the top of the spectrum of $\Sigma_G$ is studied
		in~\cite[Sec.~5.4]{cahill2025group}. The distortion also depends on the bottom of the
		spectrum, which we bound in Section~\ref{app:orderstat}.
		
		\item The smallest distortion over all choices of templates is computed
		in~\cite[Thms.~6 and~7]{mixonpacker}, where $S_d$ acts on $\mathbf 1^{\perp}$. It is
		$(\tfrac2\pi+o(1))\,d$ for $S_d$, $(\tfrac4\pi+o(1))\,d$ for $B_d$ and $D_d$, and
		$(\tfrac2\pi+o(1))\,k$ for $I_2(k)$. For the sign group it is $1$, since $a_j=e_j$ gives
		$f_A(x)=|x|$. Thus, in the limit $D\to\infty$, Gaussian templates lose a factor of order $\sqrt d$ for
		$\{\pm1\}^d$, $S_d$, $B_d$ and $D_d$, and only a factor tending to $\sqrt3$ for $I_2(k)$.
		
		\item Since $\lambda_{\max}(M_G)/d\to1$ by Theorem~\ref{thm:growth}, the upper bounds on $\kappa_G$ in
		Theorem~\ref{thm:growth} are asymptotically sharp for $S_d$ and $B_d$ if and only if
		$d^2\lambda_{\min}(M_G)$ tends to $\pi/2$, respectively $\pi/8$. The numerics of
		Section~\ref{sec:experiments} suggest that this holds, but we do not prove it.
		
		\item The chamber chart needs a group generated by reflections, so it is not available for
		multisets of $n\ge2$ points in $\R^{d'}$ with $d'\ge2$. Indeed, the group $S_n$ acting on
		$\R^{n\times d'}$ by permuting rows then contains no reflection, because the
		$(-1)$-eigenspace of every nontrivial permutation has dimension a multiple of $d'$. For
		such groups the cross-cone terms of Proposition~\ref{thm:formula} may set the lower
		constant. For max filters over reflection groups they do not (see the remark after
		Example~\ref{ex:two-piece}).
	\end{enumerate}
\end{rems}
It remains to prove Theorem~\ref{thm:growth}.

\section{Two Lemmas on Order Statistics}
\label{app:orderstat}

Let i.i.d. real random variables
$X_1,\dots,X_d$  with cdf $F$, and we write
$X_{(1)}\le\cdots\le X_{(d)}$ for their order statistics. 

Recall that $(X_{(i)})_{i=1}^d$ is a Markov chain. Note that $X_{(i)}$ has the law of $X_1$ conditioned to lie in
$(X_{(i-1)},X_{(i+1)})$, and it is conditionally independent of the other order statistics~\cite[\S2.5]{david2003order}. Here $X_{(0)}:=-\infty$ and $X_{(d+1)}:=+\infty$. Let
\[
v_i:=\operatorname{Var}\bigl(X_{(i)}\mid X_{(i-1)},X_{(i+1)}\bigr)
\quad\text{and}\quad W_i:=\EX[v_i].
\]

\begin{lem}
	\label{lem:orderstat-floor}
	Suppose that $\EX[X_1^2]<\infty$ and that the law of $X_1$ has a density bounded
	above by $c<\infty$. Then, for every $d\ge1$,
	\[
	\Cov\bigl(X_{(1)},\dots,X_{(d)}\bigr)
	\ \succeq\ \tfrac12\,\diag(W_1,\dots,W_d)
	\ \succeq\ \frac{1}{4c^{2}(d+1)(d+2)}\;I_d .
	\]
\end{lem}

\begin{proof}
	Let $f$ be the density of $X_1$. We write $\mathbf X:=(X_{(1)},\dots,X_{(d)})$, let $O$ and
	$E$ be the sets of odd and even indices in $[d]$, and let
	$\mathbf X_O:=(X_{(k)})_{k\in O}$ and $\mathbf X_E:=(X_{(k)})_{k\in E}$.

	Given $\mathbf X_O$, the coordinates $X_{(i)}$ with $i\in E$ are independent.  Each even index $i$
has odd neighbors $i\pm1$, except that $i+1$ may be the boundary index $d+1$. So,
given $\mathbf X_O$, the conditional density of $\mathbf X_E$ is proportional to
\[
\prod_{i\in E} f(x_i)\,\mathbf 1\bigl\{X_{(i-1)}<x_i<X_{(i+1)}\bigr\},
\]
and each factor depends on $x_i$ only. After normalization, this factor is
the density of $X_1$ conditioned to lie in $(X_{(i-1)},X_{(i+1)})$. Therefore the conditional variance of $X_{(i)}$ is $v_i$. Moreover, given $\mathbf X_E$, the coordinates $X_{(i)}$ with $i\in O$ are
independent, and the conditional variance of $X_{(i)}$ is $v_i$.

	 Fix $w\in\R^d$. Given $\mathbf X_O$, the sum
	$\sum_{i\in O}w_iX_{(i)}$ is fixed. Hence, by the previous paragraph,
	$\operatorname{Var}(\inner{w}{\mathbf X}\mid\mathbf X_O)=\sum_{i\in E}w_i^2v_i$, and by
	the law of total variance
	\[
	w^{\!\top}\Cov(\mathbf X)\,w=\operatorname{Var}\inner{w}{\mathbf X}
	\ \ge\ \EX\bigl[\operatorname{Var}(\inner{w}{\mathbf X}\mid \mathbf X_O)\bigr]
	=\sum_{i\in E}w_i^2\,W_i .
	\]
	Conditioning on $\mathbf X_E$ instead gives the same bound with $O$ in place of $E$.
	Averaging the two bounds gives
	$w^{\!\top}\Cov(\mathbf X)w\ge\tfrac12\sum_{i=1}^dw_i^2W_i$, which is the first
	inequality of the lemma.

	For the second inequality, it suffices to show that $W_i\ge1/(2c^2(d+1)(d+2))$ for
every $i$. Note that if a random variable $Y$ has mean $m$ and a density
bounded by $c'$, then $\mathbb P(|Y-m|<x)\le2c'x$ for every $x>0$. Hence
\[
\operatorname{Var}(Y)=\int_0^\infty2x\,\mathbb P(|Y-m|\ge x)\,dx
\ \ge\ \int_0^{1/(2c')}2x\,(1-2c'x)\,dx=\frac{1}{12c'^2}.
\]
Since $X_1$ has a density, $F$ is continuous. Therefore the variables $U_j:=F(X_j)$ are i.i.d.\
uniform on $(0,1)$, and their order statistics are $U_{(k)}=F(X_{(k)})$. We set
$U_{(0)}:=0$ and $U_{(d+1)}:=1$, and put $T_i:=U_{(i+1)}-U_{(i-1)}$. Given $X_{(i-1)}$ and
$X_{(i+1)}$, the variable $X_{(i)}$ has the density
$f(x)\,\mathbf 1\{X_{(i-1)}<x<X_{(i+1)}\}/T_i$, which is at most $c/T_i$. Therefore
$v_i\ge T_i^2/(12c^2)$. The $d+1$ spacings $U_{(k)}-U_{(k-1)}$, $1\le k\le d+1$, are exchangeable. Hence $T_i$, which is the sum of two of them, has the law of $U_{(2)}$. For $d\ge2$ this is the
$\operatorname{Beta}(2,d-1)$ law, and for $d=1$ we have $T_1=1$. In both cases
$\EX[T_i^2]=6/((d+1)(d+2))$, and
\[
W_i=\EX[v_i]\ \ge\ \frac{\EX[T_i^{2}]}{12c^{2}}=\frac{1}{2c^{2}(d+1)(d+2)} .
\]
 This gives the second inequality.
\end{proof}

For the uniform law on $(0,1)$, where $c=1$, the covariance formula
in~\cite{david2003order} gives that the inverse of the covariance matrix is
$(d+1)(d+2)\operatorname{tridiag}(-1,2,-1)$. So its smallest eigenvalue is the floor of
Lemma~\ref{lem:orderstat-floor} times $\cos^{-2}\tfrac{\pi}{2(d+1)}=1+O(d^{-2})$, and the
constant $1/4$ cannot be replaced by a larger one.

\begin{lem}
	\label{lem:central-spacing}
	Suppose that $\EX[X_1^2]<\infty$ and that the law of $X_1$ has a density $f$ with
	$f(t)\ge h\min\{F(t),1-F(t)\}$ for all $t\in\R$, for some $h>0$. Then, for every
	$d\ge2$ and every index $i$ with $d/4\le i\le3d/4$,
	\[
	\EX\bigl[(X_{(i+1)}-X_{(i)})^2\bigr]\ \le\ C_0\,d^{-2},
	\qquad C_0:=\max\bigl\{288/h^2,\ 2662\,\EX[X_1^2]\bigr\}.
	\]
\end{lem}

\begin{proof}
	 Let $U_{(k)}=F(X_{(k)})$ be the order
	statistics of the i.i.d.\ uniform variables $U_j:=F(X_j)$. We set $\Delta:=X_{(i+1)}-X_{(i)}$. If $x<y$ and $x\le t\le y$, then $F(t)\ge F(x)$ and
	$1-F(t)\ge1-F(y)$. Hence
	\[
	F(y)-F(x)=\int_x^yf(t)\,dt\ \ge\ h\,(y-x)\min\{F(x),1-F(y)\}.
	\]
	With $x=X_{(i)}$ and $y=X_{(i+1)}$ this gives $\Delta\le S/(h\min\{A,B\})$, where
	$A:=U_{(i)}$, $S:=U_{(i+1)}-U_{(i)}$ and $B:=1-U_{(i+1)}$. Therefore
	$\Delta^2\le h^{-2}S^2(A^{-2}+B^{-2})$. Notice that the vector $(A,S,B)$ has the
	$\operatorname{Dirichlet}(i,1,d-i)$ distribution, so for $i\ge3$ and $d-i\ge3$,
	\[
	\EX[S^2A^{-2}]=\frac{2}{(i-1)(i-2)},\qquad
	\EX[S^2B^{-2}]=\frac{2}{(d-i-1)(d-i-2)}.
	\]
	If $d\ge12$, then $i\ge3$ and $d-i\ge3$. Moreover $(i-1)(i-2)\ge\frac29i^2\ge d^2/72$,
	and $(d-i-1)(d-i-2)\ge d^2/72$. Hence $\EX[\Delta^2]\le288\,h^{-2}d^{-2}$. On the other hand, if $d\le11$, then
	$\Delta^2\le2X_{(i)}^2+2X_{(i+1)}^2\le2\sum_jX_j^2$, so $\EX[\Delta^2]\le2d\,\EX[X_1^2]\le2662\,\EX[X_1^2]\,d^{-2}$.
\end{proof}

\section{Proofs of the Growth Laws}\label{sec:proofs}

In this section we prove Theorem~\ref{thm:growth},
we study the sign group when $d\to\infty$ and $D/d\to\infty$, and we give a version of
Theorem~\ref{thm:spectral-limit} for a finite number of templates.

\subsection{Proof of Theorem~\ref{thm:growth}}
\label{app:bc}

\begin{proof}
	Let $\xi\sim N(0,I_d)$, $Z:=\Lambda_G(\xi)$, and let $M_G$, $\Sigma_G$ and $\mu_G$ be as
	in Definition~\ref{def:MG}, so $M_G=\Sigma_G+\mu_G\mu_G^{\!\top}$. By
	Theorem~\ref{thm:spectral-limit}, $\kappa_G^2=\lambda_{\max}(M_G)/\lambda_{\min}(M_G)$. Therefore we need upper and lower bounds on $\lambda_{\max}(M_G)$ and on $\lambda_{\min}(M_G)$.
	
	Let $G=S_d$ or $G=B_d$. Then $Z={\downarrow}\xi$ or $Z={\downarrow}|\xi|$. Thus $Z$ is the vector of order statistics, in decreasing order, of $d$ i.i.d.\ copies of a random
	variable $X$, which is standard normal for $S_d$ and half-normal for $B_d$. In both cases
	$\EX[X^2]=1$, and the density of $X$ is bounded by $c$, where $c=1/\sqrt{2\pi}$ for $S_d$ and
	$c=\sqrt{2/\pi}$ for $B_d$. Notice that the order of the coordinates does not change the
	eigenvalues of $M_G$ and $\Sigma_G$. Therefore, while we treat $S_d$ and $B_d$, we may assume that
	$Z=(X_{(1)},\dots,X_{(d)})$, with $X_{(1)}\le\cdots\le X_{(d)}$ as in
	Section~\ref{app:orderstat}.
	
	The largest eigenvalue of $M_G$ is close to $d$. Indeed,
	$\lambda_{\max}(M_G)\le\operatorname{tr}M_G=d$ by Definition~\ref{def:MG}. Moreover, $M_G\succeq\mu_G\mu_G^{\!\top}$, so
	$\lambda_{\max}(M_G)\ge\|\mu_G\|_2^2=\sum_{i=1}^d(\EX X_{(i)})^2$. By~\cite[Thm.~1]{hoeffding1953},
	$\frac1d\sum_{i=1}^d(\EX X_{(i)})^2\to\EX[X^2]=1$ as $d\to\infty$. Hence
	$\|\mu_G\|_2^2=(1-o(1))\,d$, and
	\[
	(1-o(1))\,d\ \le\ \lambda_{\max}(M_G)\ \le\ d .
	\]
	
	The smallest eigenvalue of $M_G$ is at least of order $d^{-2}$. Indeed, $\Sigma_G$ is the
	covariance matrix of $(X_{(1)},\dots,X_{(d)})$, so Lemma~\ref{lem:orderstat-floor} applies to it. Since $M_G\succeq\Sigma_G$, we get
	\begin{equation}\label{eq:lmin-lower}
		\lambda_{\min}(M_G)\ \ge\ \lambda_{\min}(\Sigma_G)\ \ge\ \frac{1}{4c^2(d+1)(d+2)}
		=\begin{cases}
			\dfrac{\pi}{2(d+1)(d+2)}, & G=S_d,\\[1.5ex]
			\dfrac{\pi}{8(d+1)(d+2)}, & G=B_d.
		\end{cases}
	\end{equation}
	
	The smallest eigenvalue of $M_G$ is also at most of order $d^{-2}$. Both populations
	satisfy the hypotheses of Lemma~\ref{lem:central-spacing} with $C_0=2662$, as shown after its
	proof. Put $C_1:=C_0$. If $d\le2$, then
	$\lambda_{\min}(M_G)\le\operatorname{tr}M_G/d=1\le4d^{-2}\le C_1d^{-2}$. If $d\ge3$, let
	$i:=\lfloor d/2\rfloor$ and $v:=\frac1{\sqrt2}(e_i-e_{i+1})$. Since $d/4\le i\le3d/4$,
	Lemma~\ref{lem:central-spacing} gives
	\[
	\lambda_{\min}(M_G)\ \le\ \inner{v}{M_Gv}=\tfrac12\EX\bigl[(X_{(i+1)}-X_{(i)})^2\bigr]
	\ \le\ C_0\,d^{-2}.
	\]
	Thus $\lambda_{\min}(M_G)\le C_1d^{-2}$ for every $d$.
	
	These bounds give $\kappa_G^2\ge(1-o(1))\,d^3/C_1$ and $\kappa_G^2\le4c^2\,d(d+1)(d+2)$.
	So $\kappa_G=\Theta(d^{3/2})$. $4c^2=2/\pi$ for $S_d$ and $4c^2=8/\pi$ for $B_d$,
	so we get the upper bounds on $\kappa_G$ in
	Theorem~\ref{thm:growth}(a) and, for $B_d$, in
	Theorem~\ref{thm:growth}(b).
	
	Now let $G=D_d$. We write
	$Y:={\downarrow}|\xi|$ and $\varsigma:=\prod_{i}\operatorname{sign}(\xi_i)$, so that
	$\Lambda_{D_d}(\xi)=(Y_1,\dots,Y_{d-1},\varsigma Y_d)$. The signs of the $\xi_i$ are
	independent of each other and of $|\xi|$, and each equals $\pm1$ with probability $1/2$. Therefore $\varsigma$ is independent of $Y$ and $\EX[\varsigma]=0$. Hence
	$\EX[Y_i\,\varsigma Y_d]=\EX[\varsigma]\,\EX[Y_iY_d]=0$ for $i<d$. Since $\varsigma^2=1$, we
	get
	\[
	M_{D_d}=\begin{pmatrix} M' & 0\\ 0 & \EX[Y_d^2]\end{pmatrix},
	\]
	where $M'$ is the leading $(d-1)\times(d-1)$ principal submatrix of
	$M_{B_d}=\EX[YY^{\!\top}]$, and $\EX[Y_d^2]$ is its last diagonal entry. Both are
	principal submatrices of $M_{B_d}$, so
	$\lambda_{\min}(M_{D_d})\ge\lambda_{\min}(M_{B_d})\ge\pi/(8(d+1)(d+2))$ by
	\eqref{eq:lmin-lower}. For the upper bound, let $d\ge3$ and $i:=\lfloor d/2\rfloor$, so
	that $i+1\le d-1$, and let $v:=\frac1{\sqrt2}(e_i-e_{i+1})$, now in the coordinates of $Y$.
	Then $\inner{v}{M_{D_d}v}=\inner{v}{M_{B_d}v}=\tfrac12\EX\bigl[(Y_i-Y_{i+1})^2\bigr]$. Since $Y$
	lists the order statistics of $|\xi_1|,\dots,|\xi_d|$ in decreasing order,
	$Y_i-Y_{i+1}=X_{(j+1)}-X_{(j)}$ with $j:=d-i$. Since $d/4\le j\le3d/4$,
	Lemma~\ref{lem:central-spacing} gives $\inner{v}{M_{D_d}v}\le C_0d^{-2}$. For $d\le2$ we use the trace as before. For the largest eigenvalue, notice
	that $\mu_{D_d}=(\EX Y_1,\dots,\EX Y_{d-1},0)$ and $(\EX Y_d)^2\le\EX[Y_d^2]\le\EX[\xi_1^2]=1$, because $Y_d=\min_k|\xi_k|\le|\xi_1|$.
	Thus $\|\mu_{D_d}\|_2^2\ge\|\mu_{B_d}\|_2^2-1=(1-o(1))\,d$, and therefore
	$(1-o(1))\,d\le\lambda_{\max}(M_{D_d})\le\operatorname{tr}M_{D_d}=d$.
	
	Finally, we will show that $S_d$ acting on $\mathbf 1^{\perp}$ has the same $\kappa_G$ as $S_d$
	acting on $\R^d$. Let $d\ge2$, $V:=\mathbf 1^{\perp}$, and $G'$ be $S_d$ acting on
	$V$, which we identify with $\R^{d-1}$. Note that $G'$ is a
	reflection group, with chamber $\{x_1\ge\cdots\ge x_d\}\cap V$ and chart
	$x\mapsto{\downarrow}x$. We write $\bar\xi:=\frac1d\sum_i\xi_i$ and
	$\xi_V:=\xi-\bar\xi\,\mathbf 1\in V$. Then $\xi_V$ is a standard Gaussian vector in $V$.
	Moreover, $\xi_V$ and $\bar\xi$ are independent.
	Since sorting commutes with adding a multiple of $\mathbf 1$, we have
	${\downarrow}\xi={\downarrow}\xi_V+\bar\xi\,\mathbf 1$, with ${\downarrow}\xi_V\in V$. Using that
	$\bar\xi$ is independent of $\xi_V$, $\EX\bar\xi=0$ and $\EX\bar\xi^{\,2}=1/d$, we obtain
	\[
	M_{S_d}=\EX\bigl[{\downarrow}\xi_V\,({\downarrow}\xi_V)^{\!\top}\bigr]+\tfrac1d\,\mathbf 1\mathbf 1^{\!\top}.
	\]
	Note that the first term is $M_{G'}$ on $V$ and vanishes on $\mathbf 1$, while the second
	term is the orthogonal projection onto $\operatorname{span}\{\mathbf 1\}$. Therefore $M_{S_d}\mathbf 1=\mathbf 1$, and $M_{S_d}$ acts on $V$ as $M_{G'}$. Furthermore,
	$\operatorname{tr}M_{G'}=\EX\|\xi_V\|_2^2=d-1=\dim V$, so
	$\lambda_{\min}(M_{G'})\le1\le\lambda_{\max}(M_{G'})$. We conclude that $M_{S_d}$ and $M_{G'}$
	have the same largest and smallest eigenvalues, and so $\kappa_{G'}=\kappa_{S_d}$.
\end{proof}
\subsection{The Sign Group when $d\to\infty$ and $D/d\to\infty$}
\label{app:sign-joint}

By Theorems~\ref{thm:spectral-limit} and~\ref{thm:exact}(a), for fixed $d$ the distortion
$\kappa(f_A)$ tends to $\kappa_G=\sqrt{1+2d/(\pi-2)}$ as $D\to\infty$. For the sign group, we now
show that $\kappa(f_A)$ has the same leading term $\sqrt{2/(\pi-2)}\,\sqrt d$, with probability
tending to one, when $d\to\infty$ and $D/d\to\infty$. Here $\Lambda_G(x)=|x|$, so the matrix $R$
of Theorem~\ref{thm:reflection-constants} has i.i.d.\ half-normal entries. The effect of a
rank-one mean on the largest singular value was studied in~\cite{silverstein1994spectral} when
$D$ and $d$ grow proportionally, and the bounds for the centered part are
standard~\cite{vershynin2018high}.

\begin{prop}
	\label{prop:sign-joint}
	Let $G=\{\pm1\}^d$, let $a_1,\dots,a_D$ be i.i.d.\ $N(0,I_d)$, and let $R$ be the
	matrix with rows $|a_j|^{\!\top}$. Its entries have mean $m_1:=\sqrt{2/\pi}$ and
	variance $\tau^2:=1-2/\pi$. When $d\to\infty$ and $D/d\to\infty$, with
	probability tending to one,
	\[
	\sigma_{\min}(R)=(1+o(1))\,\tau\sqrt{D},
	\qquad
	\sigma_{\max}(R)=(1+o(1))\,m_1\sqrt{Dd},
	\]
	so that the distortion satisfies
	$\kappa(f_A)=\bigl(1+o(1)\bigr)\sqrt{2/(\pi-2)}\,\sqrt d$.
\end{prop}

\begin{proof}
	We write $R=m_1 J+\tilde R$, where $J:=\mathbf 1_D\mathbf 1_d^{\!\top}$ and $\tilde R$
	has i.i.d.\ centered sub-Gaussian entries of variance $\tau^2$. We will show that
	$\sigma_{\max}(R)$ is close to $m_1\|J\|_{\operatorname{op}}=m_1\sqrt{Dd}$, which is the
	contribution of the rank-one mean, and that $\sigma_{\min}(R)$ is close to $\tau\sqrt D$,
	which is the scale of all singular values of $\tilde R$. For this we use two events, with
	$\varepsilon:=2C_{\mathrm V}\sqrt{d/D}$, where $C_{\mathrm V}$ is the constant in the bound
	quoted in (E1). Note that $\varepsilon\to0$, since $D/d\to\infty$.
	\begin{enumerate}
		\item[(E1)] $(1-\varepsilon)\tau\sqrt D\le\sigma_{\min}(\tilde R)\le\sigma_{\max}(\tilde R)\le(1+\varepsilon)\tau\sqrt D$.
		The rows of $\tau^{-1}\tilde R$ are i.i.d.\ and isotropic with $O(1)$
		sub-Gaussian norm. Hence all singular values of $\tilde R$ lie in
		$\tau[\sqrt D-C_{\mathrm V}(\sqrt d+t),\,\sqrt D+C_{\mathrm V}(\sqrt d+t)]$ with probability at least
		$1-2e^{-t^2}$~\cite[Thm.~4.6.1]{vershynin2018high}. We take $t=\sqrt d$.
		\item[(E2)] $\|s\|_2^2\le\tau^2D\sqrt{Dd}$, where $s:=\tilde R^{\!\top}\mathbf 1_D$.
		Each coordinate $s_i=\sum_j\tilde R_{ji}$ is a sum of $D$ independent centered
		variables of variance $\tau^2$, thus $\EX\|s\|_2^2=\tau^2Dd$. By Markov's inequality,
		(E2) fails with probability at most $\sqrt{d/D}$.
	\end{enumerate}
By the union bound, (E1) and (E2) both hold with probability at least
$1-2e^{-d}-\sqrt{d/D}$, and this tends to one. 
	By the triangle inequality and (E1), we have
	$\sigma_{\max}(R)\le m_1\|J\|_{\operatorname{op}}+\sigma_{\max}(\tilde R)=m_1\sqrt{Dd}+O(\tau\sqrt D)$,
	and
	\[
	\sigma_{\max}(R)\ \ge\ \frac{\|R\mathbf 1_d\|_2}{\sqrt d}
	=\Bigl\|m_1\sqrt d\,\mathbf 1_D+\tfrac1{\sqrt d}\tilde R\mathbf 1_d\Bigr\|_2
	\ \ge\ m_1\sqrt{Dd}-\sigma_{\max}(\tilde R).
	\]
	Since $\tau\sqrt D=o(m_1\sqrt{Dd})$, both bounds give
	$\sigma_{\max}(R)=(1+o(1))m_1\sqrt{Dd}$.

	For a unit vector $w\perp\mathbf 1_d$ we have $Jw=0$. Thus $Rw=\tilde Rw$, and by (E1) we
	get $\sigma_{\min}(R)\le(1+\varepsilon)\tau\sqrt D$. For the lower bound, let $\Pi:=\frac1D\mathbf 1_D\mathbf 1_D^{\!\top}$ be the
	orthogonal projector onto $\operatorname{span}\{\mathbf 1_D\}$. Notice that
	$\operatorname{span}\{\mathbf 1_D\}$ is the range of $J$. Therefore $(I_D-\Pi)J=0$, and hence
	$(I_D-\Pi)R=(I_D-\Pi)\tilde R$. Consequently, for every unit $v\in\R^d$,
	\[
	\begin{aligned}
	\|Rv\|_2\ &\ge\ \|(I_D-\Pi)\tilde Rv\|_2
	\ \ge\ \|\tilde Rv\|_2-\|\Pi\tilde Rv\|_2\\
	&=\|\tilde Rv\|_2-\frac{|\inner{s}{v}|}{\sqrt D}
	\ \ge\ \sigma_{\min}(\tilde R)-\frac{\|s\|_2}{\sqrt D},
	\end{aligned}
	\]
	where we used $\|\Pi\tilde Rv\|_2=|\mathbf 1_D^{\!\top}\tilde Rv|/\sqrt D$ and the
	Cauchy Schwarz inequality. By (E1) and (E2), the right-hand side is at least
	$(1-\varepsilon)\tau\sqrt D-\tau(Dd)^{1/4}=\bigl(1-\varepsilon-(d/D)^{1/4}\bigr)\tau\sqrt D$,
	uniformly in $v$. Since $\varepsilon\to0$ and $d/D\to0$, this is $(1-o(1))\tau\sqrt D$.
	Thus $\sigma_{\min}(R)=(1+o(1))\tau\sqrt D$.

	In particular $\sigma_{\min}(R)>0$, so $\operatorname{rank}R=d$. Now the distortion
	statement follows from Theorem~\ref{thm:reflection-constants} and
	$m_1/\tau=\sqrt{2/(\pi-2)}$.
\end{proof}

\subsection{A Finite Number of Templates}
\label{app:finiteD}

\begin{prop}
	\label{prop:finiteD}
	Let $G\le O(d)$ be a finite reflection group, let $\varepsilon,p\in(0,1)$, let
	$a_1,\dots,a_D$ be i.i.d.\ $N(0,I_d)$, and $R$ be the matrix of Theorem~\ref{thm:reflection-constants} for these
	templates. If
	\[
	D\ \ge\ 1000\,\varepsilon^{-2}\,d\,\log^2\!\Bigl(\frac{1100\,d}{\varepsilon p}\Bigr)
	\log\Bigl(\frac{3d}{p}\Bigr),
	\]
	then with probability at least $1-p$,
	\[
	(1-\varepsilon)\,D\,M_G\ \preceq\ R^{\!\top}R\ \preceq\ (1+\varepsilon)\,D\,M_G .
	\]
	In that case $f_A$ is bi-Lipschitz, $a(f_A)^2$ and $b(f_A)^2$ lie within the factors
	$1\pm\varepsilon$ of $D\lambda_{\min}(M_G)$ and $D\lambda_{\max}(M_G)$, and
	\[
	\Bigl(\frac{1-\varepsilon}{1+\varepsilon}\Bigr)^{1/2}\kappa_G\ \le\ \kappa(f_A)\ \le\
	\Bigl(\frac{1+\varepsilon}{1-\varepsilon}\Bigr)^{1/2}\kappa_G .
	\]
\end{prop}

 For $d\ge2$ and $\varepsilon=1/2$ we get $\kappa_G/\sqrt3\le\kappa(f_A)\le\sqrt3\,\kappa_G$
with probability at least $1-p$, if $D\ge C_\ast\,d\log^3(d/p)$ for a constant
$C_\ast$. Note also that this condition is independent on $G$. On the other hand, if $D<d$, then
$\operatorname{rank}R<d$, and by Theorem~\ref{thm:reflection-constants} the map $f_A$ does not
separate orbits. Therefore the power of $d$ is optimal.

\begin{proof}
	Let $M:=M_G$ and $Z_j:=\Lambda_G(a_j)$, so that $R^{\!\top}R=\sum_jZ_jZ_j^{\!\top}$. By
	Theorem~\ref{thm:spectral-limit} we have $M\succ0$. Let $W_j:=M^{-1/2}Z_j$ and
	$S:=\frac1D\sum_jW_jW_j^{\!\top}$. Then $\EX S=I$ and $R^{\!\top}R/D=M^{1/2}SM^{1/2}$. Hence it
	suffices to show that $(1-\varepsilon)I\preceq S\preceq(1+\varepsilon)I$ with probability
	at least $1-p$.
	
	Let $\xi\sim N(0,I_d)$ and
	$Z:=\Lambda_G(\xi)$, so that each $Z_j$ has the law of $Z$. Let $B\subseteq\operatorname{int}C$
	be a Borel set. Then $Z\in B$ if and only if $\xi\in U_gB$ for some $g\in G$. Notice that the
	sets $U_gB$, $g\in G$, are pairwise disjoint, since each orbit meets $C$ in exactly one point
	and every point of $\operatorname{int}C$ has trivial stabilizer~\cite[Thm.~1.12]{humphreys}.
	The standard Gaussian density $\phi_d$ is invariant under $O(d)$, therefore
	$\mathbb{P}(Z\in B)=\sum_{g\in G}\mathbb{P}(\xi\in U_gB)=|G|\,\mathbb{P}(\xi\in B)$. Moreover,
	the boundary of $C$ lies in finitely many hyperplanes, so it has Lebesgue measure zero. Since
	$Z\in\partial C$ only if $\xi\in U_g(\partial C)$ for some $g\in G$, we have
	$\mathbb{P}(Z\in\partial C)\le|G|\,\mathbb{P}(\xi\in\partial C)=0$. As $Z\in C$, we
	conclude that $Z$ has density $|G|\,\phi_d\mathbf 1_C$. This density is log-concave, because
	$\phi_d$ is log-concave and $C$ is convex. Consequently $W:=M^{-1/2}Z$, being an invertible
	linear image of $Z$, also has a log-concave density. Furthermore
	$\EX WW^{\!\top}=M^{-1/2}MM^{-1/2}=I$, so $\EX\|W\|_2^2=d$.
	
	Now let $\mathcal B:=\{x:\|x\|_2\le3\sqrt d\}$. By Markov's
	inequality we have $\theta:=\mathbb{P}(W\in\mathcal B)\ge8/9$. Since the
	law of $W$ has a log-concave density, it is a log-concave measure~\cite{borell1975convex}.
	Therefore Borell's lemma~\cite{borell1974convex} for the symmetric convex set $\mathcal B$ gives
	$\mathbb{P}(W\notin t\mathcal B)\le\theta\bigl((1-\theta)/\theta\bigr)^{(t+1)/2}
	=\theta^{(1-t)/2}(1-\theta)^{(t+1)/2}$ for $t\ge1$. The right-hand side is nonincreasing in
	$\theta$, so for $t\ge1$,
	\[
	\mathbb{P}\bigl(\|W\|_2>3t\sqrt d\bigr)\ \le\ \tfrac89\bigl(\tfrac18\bigr)^{(t+1)/2}\ \le\ e^{-t}.
	\]
	We write $\zeta:=\|W\|_2/(3\sqrt d)$. Then for $t\ge1$ we have
	\[
	\EX\bigl[\zeta^2\mathbf 1\{\zeta>t\}\bigr]=t^2\,\mathbb{P}(\zeta>t)+\int_t^\infty2s\,\mathbb{P}(\zeta>s)\,ds
	\ \le\ (t^2+2t+2)\,e^{-t}\ \le\ 5t^2e^{-t},
	\]
	since $4t^2-2t-2=2(2t+1)(t-1)\ge0$. Multiplying by $9d$, we obtain
	\begin{equation}\label{eq:W-tail}
		\EX\bigl[\|W\|_2^2\mathbf 1\{\|W\|_2>3t\sqrt d\}\bigr]\ \le\ 45\,t^2d\,e^{-t}
		\qquad\text{for }t\ge1.
	\end{equation}
	
	 Let $N:=1100\,d/(\varepsilon p)$, fix $t:=2\log N\ge1$, and let
	$E_j:=\{\|W_j\|_2\le3t\sqrt d\}$. On $E_j$ we have $\|W_j\|_2^2\le9\,t^2d$. Moreover,
	$e^{-t}=N^{-2}$ and $\log^2N\le N$, so by~\eqref{eq:W-tail}
	\[
	\eta:=\EX\bigl[\|W_j\|_2^2\mathbf 1_{E_j^c}\bigr]\ \le\ \frac{180\,d\log^2N}{N^2}
	\ \le\ \frac{180\,d}{N}\ =\ \frac{180}{1100}\,\varepsilon p\ \le\ \frac{\varepsilon p}{6}.
	\]
	We denote $\Delta:=\EX[W_jW_j^{\!\top}\mathbf 1_{E_j^c}]$.
	Notice that $0\preceq\Delta\preceq\eta I$. We also write
	\[
	S':=\frac1D\sum_jW_jW_j^{\!\top}\mathbf 1_{E_j},\qquad
	S'':=S'+\Delta=\frac1D\sum_j\bigl(W_jW_j^{\!\top}\mathbf 1_{E_j}+\Delta\bigr).
	\]
	Then $\EX S''=I$. Furthermore, the $D$ summands of $S''$ are independent and positive
	semidefinite, and their largest eigenvalue is at most $L:=(9t^2d+\eta)/D\le10\,t^2d/D$, since
	$\eta<1\le t^2d$.
	
	We apply~\cite[Thm.~1.1]{tropp2012user} with $\delta:=\varepsilon/2$, together with
	the elementary bounds $e^{-\delta}(1-\delta)^{\delta-1}\le e^{-\delta^2/2}$ and
	$e^{\delta}(1+\delta)^{-1-\delta}\le e^{-\delta^2/3}$ for $\delta\in[0,1]$. We obtain
	\[
	\mathbb{P}\bigl(\lambda_{\min}(S'')\le1-\delta\bigr)
	+\mathbb{P}\bigl(\lambda_{\max}(S'')\ge1+\delta\bigr)\ \le\ 2d\,e^{-\delta^2/(3L)}
	\ \le\ \frac{2p}{3}.
	\]
	For the last inequality, notice that $L\le10t^2d/D=40\,d\log^2N/D$ and $\delta=\varepsilon/2$
	give $\delta^2/(3L)\ge\varepsilon^2D/(480\,d\log^2N)\ge\log(3d/p)$.
	On the other hand, the random variable
	$T:=\frac1D\sum_j\|W_j\|_2^2\mathbf 1_{E_j^c}$ has $\EX T=\eta$, so by Markov's inequality
	$\mathbb{P}(T\ge\varepsilon/2)\le2\eta/\varepsilon\le p/3$.
	
	Therefore, we get that with probability at least $1-p$ we have
	$(1-\frac\varepsilon2)I\preceq S''\preceq(1+\frac\varepsilon2)I$ and $T<\varepsilon/2$. Then
	$S\succeq S'=S''-\Delta\succeq(1-\frac\varepsilon2-\eta)I\succeq(1-\varepsilon)I$. Moreover,
	$S-S'\succeq0$ has trace $T$ and $S'\preceq S''$, therefore
	$\lambda_{\max}(S)\le\lambda_{\max}(S')+T\le\lambda_{\max}(S'')+T<1+\varepsilon$. Thus
	$(1-\varepsilon)I\preceq S\preceq(1+\varepsilon)I$, and consequently
	$(1-\varepsilon)DM\preceq R^{\!\top}R\preceq(1+\varepsilon)DM$. Moreover,
	$\lambda_{\min}(R^{\!\top}R)>0$, so $\operatorname{rank}R=d$. We get the final result
	from Theorem~\ref{thm:reflection-constants}, since $a(f_A)^2=\lambda_{\min}(R^{\!\top}R)$ and
	$b(f_A)^2=\lambda_{\max}(R^{\!\top}R)$.
\end{proof}

\section{Optimal Constants of Piecewise-Linear Invariants}\label{sec:formula}

In this section we study the optimal constants $a(f)$ and $b(f)$ of the maps in our
class, for any finite orthogonal group $G$. 
\subsection{Structure of the Maps}
\label{sec:homogeneous-pl}

	Throughout, $G$ is a finite group acting orthogonally on $\R^d$ via $U_g$, and $f:\R^d\to\R^m$ is as in Definition~\ref{piecelinear}. By the remark after Definition~\ref{piecelinear}, its linear components $f_1,\ldots,f_N:\R^d\to\R^m$ satisfy
	\[
	f(x)\in\{f_1(x),\ldots,f_N(x)\} \qquad \forall x\in\R^d.
	\]
	We denote the maximum operator norm of these pieces by
	\[
	L_f:=\max_{i\in[N]}\|f_i\|_{\operatorname{op}}.
	\]

	\begin{prop}\label{prop:conic-decomp}
	Let $f:\R^d\to\R^m$ be continuous, and let $f_1,\ldots,f_N:\R^d\to\R^m$ be linear
	maps with $f(x)\in\{f_1(x),\ldots,f_N(x)\}$ for every $x\in\R^d$. Then there is a
	finite collection $\mathcal H$ of linear hyperplanes in $\R^d$ with the following
	property. Every connected component of $\R^d\setminus\bigcup_{H\in\mathcal H}H$ is an
	open convex polyhedral cone. We write $\{C_\alpha\}_{\alpha\in I}$ for these cones,
	where $I$ is a finite index set. On each $C_\alpha$, $f$ coincides with one of the
	maps $f_i$. In particular, every point outside the union of the hyperplanes in
	$\mathcal H$ has a neighborhood on which $f$ agrees with one of the maps $f_i$.
\end{prop}

\begin{proof}Without loss of generality, we may assume that $f_1,\ldots,f_N$ are pairwise distinct. For
$i<j$ and $k\in[m]$, let $\ell_{ijk}$ be the $k$-th coordinate of $f_i-f_j$,
a linear functional on $\R^d$. Let $\mathcal H$ be the set of kernels of the nonzero
functionals $\ell_{ijk}$. Each of these kernels is a hyperplane of dimension $d-1$.
Each connected component of $\R^d\setminus\bigcup_{H\in\mathcal H}H$ is an
open convex polyhedral cone. There are finitely many components. We fix one of them, $\Omega$.

Let $z\in \Omega$ and $i<j$. Since $f_i\neq f_j$, some $\ell_{ijk}$ is nonzero. Its kernel
belongs to $\mathcal H$, so it does not contain $z$. Hence $\ell_{ijk}(z)\neq0$, and
$f_i(z)\neq f_j(z)$.

For $i\in[N]$ we put $S_i:=\{z\in \Omega: f(z)=f_i(z)\}$. These sets cover $\Omega$, because $f(z)$ is always one
of the values $f_i(z)$. They are closed in $\Omega$, because $f$ and $f_i$ are
continuous, and they are pairwise disjoint by the previous paragraph. Therefore each $S_i$ is also open in
$\Omega$, as the complement of the union of the others. Since $\Omega$ is convex, it is
connected. Hence exactly one $S_i$ is nonempty, and it equals $\Omega$.
\end{proof}

	Since $f$ is $G$-invariant, it descends to a well-defined map
	$\widehat f:\R^d/G\to\R^m$, given by $\widehat f([x])=f(x)$. 
	\begin{prop}
		\label{prop:upper-lipschitz}
		For all $x,y\in\R^d$,
		\[
		\|\widehat f([x])-\widehat f([y])\|_2 \leq L_f\,d_G([x],[y]).
		\]
	\end{prop}

	\begin{proof}
	Fix $x,y\in\R^d$ and
		put $\gamma(t):=(1-t)x+ty$ for $t\in[0,1]$. The union of the hyperplanes in
		$\mathcal H$ has empty interior, so the union of the open cones $C_\alpha$ of
		Proposition~\ref{prop:conic-decomp} is dense in $\R^d$. Since there are finitely many
		cones, their closures $\overline{C_\alpha}$ cover $\R^d$, and by continuity $f$ agrees with a single
		linear component on each of them. Since each $\overline{C_\alpha}$ is closed and
		convex, the set $J_\alpha:=\{t\in[0,1]:\gamma(t)\in\overline{C_\alpha}\}$ is empty, a
		single point, or a closed interval. These finitely many sets
		cover $[0,1]$. Let $0=t_0<t_1<\dots<t_r=1$ be the points $0$, $1$ and the
		endpoints of the nonempty sets $J_\alpha$ that lie in $(0,1)$. For each $\ell$, the
		midpoint of $[t_{\ell-1},t_\ell]$ lies in some $J_\alpha$, and no endpoint of
		$J_\alpha$ lies strictly between $t_{\ell-1}$ and $t_\ell$. Therefore
		$[t_{\ell-1},t_\ell]\subseteq J_\alpha$, and $f$ agrees on $\gamma([t_{\ell-1},t_\ell])$ with
		the linear component $f_{i_\ell}$ that it agrees with on $\overline{C_\alpha}$. Hence
		\[
		\|f(\gamma(t_\ell))-f(\gamma(t_{\ell-1}))\|_2
		=\|f_{i_\ell}\bigl(\gamma(t_\ell)-\gamma(t_{\ell-1})\bigr)\|_2
		\le L_f(t_\ell-t_{\ell-1})\|x-y\|_2 .
		\]
		Summing over $\ell$, we obtain $\|f(x)-f(y)\|_2\le L_f\|x-y\|_2$. Notice that, since we
		use the closed cones, the argument also covers a segment that lies inside a hyperplane of
		$\mathcal H$.
		Since $G$ is finite, we may choose
		$g\in G$ with $d_G([x],[y])=\|x-U_gy\|_2$. By $G$-invariance of $f$, we conclude that
		\[
		\|\widehat f([x])-\widehat f([y])\|_2=\|f(x)-f(U_gy)\|_2\le L_f\|x-U_gy\|_2
		=L_f\,d_G([x],[y]).
		\]
	\end{proof}

\begin{rem*}[Bias-free ReLU networks]
	Let $\sigma(t):=\max\{t,0\}$ act componentwise, and
	$\mathcal N(x):=W_L\sigma\bigl(W_{L-1}\sigma(\cdots\sigma(W_1x)\cdots)\bigr)$ be a
	ReLU network without biases. Its $G$ symmetrization
	$f(x):=\frac1{|G|}\sum_{g\in G}\mathcal N(U_gx)$ lies in our class. Indeed, $f$ is continuous
	and $G$-invariant. Moreover, it is positively homogeneous, since $\sigma(\lambda t)=\lambda\sigma(t)$
	for $\lambda\ge0$. On each region where the activation patterns of the hidden units of all
	the networks $x\mapsto\mathcal N(U_gx)$, $g\in G$, are fixed, $f$ is an average of linear maps,
	hence linear. There are finitely many such regions. Therefore, if $f$ is injective on
	$\R^d/G$, then by Theorem~\ref{thm:main} it is bi-Lipschitz. The same argument applies to any
	bias-free network with ReLU, max or sorting layers that is $G$-invariant by construction.

	Injectivity depends on the pooling. For example, let $G=S_d$ with
$d\ge2$, and let $\mathcal N$ be linear. Averaging gives $f=\mathcal N\Pi$, where $\Pi$ is
the orthogonal projection onto $\operatorname{span}\{\mathbf 1\}$. So $f$ depends only on the
mean of the coordinates of $x$, and it is not injective. Max-pooling over $G$ instead gives the
max filter of Section~\ref{sec:random-maxfilter}, which is injective for enough generic
templates~\cite{cahill2025group}. Finally, with biases, injectivity no longer implies a lower
Lipschitz bound~\cite[Remark~6]{balan2024sorting}.
\end{rem*}

\begin{rem*}[Sorting-based embeddings]
	Sorting-based embeddings defined in~\cite{balan2024sorting} belong to our class, since sorting the
	coorbit $(\inner{x}{U_g\phi})_{g\in G}$ of a template $\phi$ is continuous, positively
	homogeneous, piecewise linear and $G$-invariant. On the other hand, our class is
	strictly larger than the family of sorting architectures. To see this, let $S_2$ act
	on $\R^2$ by permuting the coordinates. For a template $\phi=(a,b)$, the two sorted entries
	of the coorbit $(ax+by,\,bx+ay)$ are $\frac12\bigl((a+b)(x+y)\pm|a-b|\,|x-y|\bigr)$.
	Every component of maps built from finitely many sorting blocks, where each block sorts
	the coorbits of finitely many templates separately, lies in
	$\mathcal W:=\operatorname{span}\{x+y,\,|x-y|\}$. Now the map
	$(x,y)\mapsto\bigl(x+y,\max(|x|,|y|)\bigr)$ lies in our class, and it is injective on
	$\R^2/S_2$. Indeed, for $x\ge y$, with $s:=x+y$ and $m:=\max(|x|,|y|)$, we have
	$x=m-\max\{-s,0\}$ and $y=s-x$. However, $\max(|x|,|y|)\notin\mathcal W$, because an
	identity $\max(|x|,|y|)=\alpha(x+y)+\beta|x-y|$ at the points $(1,0)$, $(-1,0)$ and $(1,1)$
	gives $\alpha+\beta=1$, $-\alpha+\beta=1$ and $2\alpha=1$, which is contradiction.
\end{rem*}

\subsection{The Lower Constant as a Minimum over Pairs of Cones}

The constant $a(f)$ is an infimum over pairs of orbits. In this subsection we study it as finitely many local problems, one for each pair of cones. A cone-by-cone lower bound of
this kind was obtained for max filtering in~\cite[Thm.~26, Cor.~30]{mixonqaddura}.

 If $x\in C_{\alpha}$ and $y\in C_{\beta}$ with
$\alpha\neq\beta$, then $f(x)-f(y)=A_{\alpha}x-A_{\beta}y$ involves two different matrices, and it
can be small even when $A_{\alpha}$ and $A_{\beta}$ are both well conditioned. So $a(f)$ can be
smaller than $\min_{\alpha}\sigma_{\min}(A_{\alpha})$, as Example~\ref{ex:two-piece} shows.
Therefore Definition~\ref{def:mu} uses all pairs $(\alpha,\beta)$, with two constraints:
$x\in C_{\alpha}$ and $y\in C_{\beta}$, and $y$ is a nearest point of its own orbit to $x$. These
constraints make the resulting minimum equal to $a(f)$ (Proposition~\ref{thm:formula}). For max
filters over reflection groups, the pairs of distinct cones do not lower the minimum (see the
remark after Example~\ref{ex:two-piece}).

\begin{defin}\label{def:mu}
	For $\alpha,\beta\in I$, the \emph{cross-cone modulus} is
	\[
	\begin{aligned}
	\nu(\alpha,\beta):=\inf\biggl\{\ &\frac{\|A_{\alpha}x-A_{\beta}y\|_2}{\|x-y\|_2}\\
	&\;:\; x\in C_{\alpha},\ y\in C_{\beta},\ x\neq y,\ \|x-y\|_2=d_{G}(x,y)\biggr\}.
	\end{aligned}
	\]
	\end{defin}

Since the action is orthogonal, the constraint $\|x-y\|_2=d_{G}(x,y)$ is equivalent to
$\inner{x}{U_{g}y-y}\le0$ for every $g\in G$. The feasible set is never empty, because
$x=0$ and any $y\neq0$ in $C_{\beta}$ satisfy these constraints. Therefore $\nu(\alpha,\beta)$ is finite. We show in Proposition~\ref{prop:mu-attained} that the infimum is attained.

\begin{prop}\label{thm:formula}
	\leavevmode
	\begin{enumerate}
		\item[(i)] $a(f)=\min_{\alpha,\beta\in I}\nu(\alpha,\beta)$, the minimum over
		the finite set $I\times I$.
		\item[(ii)] For every $\alpha\in I$,
		$\nu(\alpha,\alpha)=\sigma_{\min}(A_{\alpha})$, and this infimum is a minimum.
		\item[(iii)] $b(f)=\max_{\alpha\in I}\sigma_{\max}(A_{\alpha})$.
	\end{enumerate}
\end{prop}

\begin{proof}
	\emph{(i)} Every $\nu(\alpha,\beta)$ is finite, as noted after Definition~\ref{def:mu},
	and $I$ is finite. Therefore $a:=\min_{\alpha,\beta}\nu(\alpha,\beta)$ exists and $a<\infty$.

	First we will show that $a\ge a(f)$. Let $(x,y)$ be feasible for $\nu(\alpha,\beta)$. Then $x\neq y$
	and $\|x-y\|_2=d_{G}(x,y)>0$, so $[x]\neq[y]$. Since $f$ agrees with $A_{\alpha}$ on
	$C_{\alpha}$ and with $A_{\beta}$ on $C_{\beta}$,
	\[
	\frac{\|A_{\alpha}x-A_{\beta}y\|_2}{\|x-y\|_2}=\frac{\|f(x)-f(y)\|_2}{d_{G}(x,y)}\ \ge\ a(f).
	\]
	Taking the infimum gives $\nu(\alpha,\beta)\ge a(f)$ for every pair, so $a\ge a(f)$.

	Conversely, let $[x]\neq[y]$. We choose $g\in G$ with $d_{G}(x,y)=\|x-U_{g}y\|_2$, and we put
	$y':=U_{g}y$. Then $f(y')=f(y)$ and $\|x-y'\|_2=d_{G}(x,y')=d_{G}(x,y)>0$. Since  $C_\gamma$ cover $\R^d$, there are $\alpha,\beta\in I$ with
	$x\in C_{\alpha}$ and $y'\in C_{\beta}$. Notice that the pair $(x,y')$ is feasible for
	$\nu(\alpha,\beta)$, and therefore
	\[
	\frac{\|f(x)-f(y)\|_2}{d_{G}(x,y)}=\frac{\|A_{\alpha}x-A_{\beta}y'\|_2}{\|x-y'\|_2}
	\ \ge\ \nu(\alpha,\beta)\ \ge\ a.
	\]
	Taking the infimum over $[x]\neq[y]$ gives $a(f)\ge a$. Hence $a(f)=a$, which proves (i).

	\emph{(ii)} For $x\neq y$ in $C_{\alpha}$ we have
	$\|A_{\alpha}(x-y)\|_2/\|x-y\|_2\ge\sigma_{\min}(A_{\alpha})$, so
	$\nu(\alpha,\alpha)\ge\sigma_{\min}(A_{\alpha})$. For the reverse inequality, let $v$ be a
	unit vector with $\|A_{\alpha}v\|_2=\sigma_{\min}(A_{\alpha})$, let $p$ be an interior
	point of $C_{\alpha}$, and put $x:=p+\varepsilon v$ and $y:=p$, with $\varepsilon>0$ so
	small that $x\in C_{\alpha}$. If the orbit of $p$ has points other than $p$, we also take $\varepsilon$ less than half the distance from $p$ to these points. Then no point of the orbit of $y$ is closer to $x$ than
	$y$ itself. Thus $\|x-y\|_2=d_{G}(x,y)$, and the pair is feasible. The ratio at this pair is
	$\|A_{\alpha}v\|_2=\sigma_{\min}(A_{\alpha})$, so the infimum is attended.

	\emph{(iii)} The first part of the proof of Proposition~\ref{prop:upper-lipschitz},
	which shows that $f$ is $L_f$-Lipschitz on $\R^d$, only uses the linear components that
	$f$ agrees with on the closed cones $C_{\gamma}$, and these are the maps $A_{\gamma}$. Therefore the same proof, with $L_f$ replaced by
	$\max_{\gamma\in I}\|A_{\gamma}\|_{\operatorname{op}}=\max_{\gamma\in I}\sigma_{\max}(A_{\gamma})$,
	gives $b(f)\le\max_{\gamma\in I}\sigma_{\max}(A_{\gamma})\leq L_f$. This bound is an equality. Indeed, choose $\alpha\in I$ with
	$\sigma_{\max}(A_\alpha)=\max_{\gamma\in I}\sigma_{\max}(A_\gamma)$, let $p$ be an interior
	point of $C_{\alpha}$, and put $x=p+\varepsilon v$ and $y=p$, where $v$ is a unit vector
	with $\|A_{\alpha}v\|_2=\sigma_{\max}(A_{\alpha})$ and $\varepsilon>0$ is chosen as in (ii).
	Then $\|f(x)-f(y)\|_2/d_{G}(x,y)=\|A_{\alpha}v\|_2=\sigma_{\max}(A_{\alpha})$, so
	$b(f)\ge\sigma_{\max}(A_{\alpha})$. This equality
	was sketched for general maps of this kind in a talk~\cite{ghoreishi2024maxfilter}, as reported
	in~\cite[Sec.~3.1]{mixonqaddura}. For max filters it is~\cite[Thm.~14]{mixonqaddura}.
\end{proof}

\begin{ex}\label{ex:two-piece}
	Let $G$ be trivial, acting on $\R$, let $h>0$, and consider $f_h\colon\R\to\R^2$ given by
	$f_h(t):=(t,0)$ for $t\ge0$ and $f_h(t):=(-t,ht)$ for $t\le0$. Notice that $f_h$ is in our
	class and that it is injective. Its cones are $[0,\infty)$ and $(-\infty,0]$, and since $G$ is
	trivial, the constraint $\|x-y\|_2=d_G(x,y)$ in Definition~\ref{def:mu} is always satisfied.
	The smallest singular values of the two pieces are $1$ and $\sqrt{1+h^2}$. Therefore, by
	Proposition~\ref{thm:formula}(ii),(iii), the two same-cone moduli are $1$ and
	$\sqrt{1+h^2}$, and $b(f_h)=\sqrt{1+h^2}$. For the cross-cone moduli we take $x\ge0$ and
	$y=-s\le0$. We may assume that
	$x+s=1$. Then $|x-y|=1$ and
	\[
	\|f_h(x)-f_h(y)\|_2^2=(x-s)^2+h^2s^2=(1-2s)^2+h^2s^2 .
	\]
	The right hand side is minimal at $s=2/(4+h^2)\in(0,1)$, where it equals $h^2/(4+h^2)$.
	If we exchange  $x$ and $y$ we get the same value. Thus both cross cone moduli
	are equal to $h/\sqrt{4+h^2}$, and from Proposition~\ref{thm:formula}(i) we obtain
	\[
	a(f_h)=\frac{h}{\sqrt{4+h^2}},\qquad
	\kappa(f_h)=\frac{\sqrt{(1+h^2)(4+h^2)}}{h}.
	\]
	Consequently, the lower constant is determined by the cross-cone terms. Moreover, the
	distortion tends to infinity as $h\to0$, and in this limit $f_h$ tends to the map
	$t\mapsto(|t|,0)$, which is not injective.
\end{ex}

\begin{rem*}
	In finite reflection groups, the cross-cone
	moduli do not change the optimal lower Lipschitz constant.
	Indeed, let $f_A=R\Lambda_G$ as in
	Theorem~\ref{thm:reflection-constants}, and let $C$ be the
	fundamental chamber. Since the open cones of
	Proposition~\ref{prop:conic-decomp} have dense union, at least
	one of them intersects the interior of $C$. Let $C_\alpha$
	be its closure. On this intersection, $\Lambda_G(x)=x$,
	so $A_\alpha x=f_A(x)=Rx$. Since the intersection is open,
	we have $A_\alpha=R$. Therefore, by
	Proposition~\ref{thm:formula}(ii) and
	Theorem~\ref{thm:reflection-constants},
	\[
	\nu(\alpha,\alpha)=\sigma_{\min}(R)=a(f_A).
	\]
	So, the optimal lower constant is attained within a
	single cone. This was already known for max filters
	from~\cite[Lem.~9]{mixonpacker}, and is also studied in~\cite[Lem.~27]{mixonqaddura}.
\end{rem*}

\subsection{Existence of a Minimizing Pair}

\begin{prop}\label{prop:mu-attained}
	Let $\alpha,\beta\in I$. Then the infimum defining $\nu(\alpha,\beta)$ is a minimum.
\end{prop}

\begin{proof}
	We prove the following more general statement.
	
	\noindent\emph{Claim.}
	Let $V\subseteq\R^d$ be a linear subspace, let
	$K_1\subseteq\R^d$ and $K_2\subseteq V$ be closed convex
	polyhedral cones, let $\Gamma$ be a finite group acting
	orthogonally on $\R^d$, and let
	$B_1,B_2:\R^d\to\R^m$ be linear. Let
	\[
	\mathcal F:=\bigl\{(u,v)\in K_1\times K_2:
	u\neq v,\ 
	\inner{u}{U_gv-v}\le0\ \text{for all }g\in\Gamma\bigr\},
	\]
	and let
	\[
	\rho(u,v):=\frac{\|B_1u-B_2v\|_2}{\|u-v\|_2}.
	\]
	If $\mathcal F\neq\emptyset$, then
	$\inf_{\mathcal F}\rho$ is a minimum.
	
	The proposition follows by taking $V=\R^d$,
	$K_1=C_\alpha$, $K_2=C_\beta$,
	$\Gamma=G$, $B_1=A_\alpha$ and $B_2=A_\beta$.
	Indeed, as noted after Definition~\ref{def:mu},
	the constraint $\|x-y\|_2=d_G(x,y)$ is equivalent to
	$\inner{x}{U_gy-y}\le0$ for all $g\in G$,
	and the feasible set is nonempty.
	
	We prove the Claim by induction on $\dim V$.
	Let $\nu:=\inf_{\mathcal F}\rho$.
	Both $\mathcal F$ and $\rho$ are invariant under
	the scaling $(u,v)\mapsto(tu,tv)$, $t>0$.
	Hence we can choose a minimizing sequence
	$(u_n,v_n)\in\mathcal F$ with
	\[
	\|u_n\|_2^2+\|v_n\|_2^2=1,
	\]
	and pass to a subsequence such that
	$(u_n,v_n)\to(u_*,v_*)$.
	If $u_*\neq v_*$, then $(u_*,v_*)\in\mathcal F$
	and $\rho(u_*,v_*)=\nu$ by continuity.
	Thus the infimum is attained.
	If $\dim V=0$, then $v_n=0$ and
	$\|u_*\|_2=1$, so $u_*\neq v_*$.
	This proves the Claim when $\dim V=0$.
	
	Now assume the induction hypothesis and that $u_*=v_*=:p$.
	Then $\|p\|_2=1/\sqrt2$ and $p\in K_1\cap K_2$.
	We first show that $B_1p=B_2p$.
	Otherwise,
	\[
	\|B_1u_n-B_2v_n\|_2
	\longrightarrow\|(B_1-B_2)p\|_2>0,
	\]
	while $\|u_n-v_n\|_2\to0$.
	This gives $\rho(u_n,v_n)\to\infty$,
	contrary to $\rho(u_n,v_n)\to\nu<\infty$.
	Consequently,
	\[
	\rho(p+u,p+v)=\rho(u,v)
	\qquad (u\neq v).
	\]
	
	Write
	\[
	K_1=\{u\in\R^d:
	\inner{a_i}{u}\le0,\ i=1,\dots,k\},
	\]
	and let $T_1$ be the cone defined by the
	constraints with $\inner{a_i}{p}=0$.
	The remaining constraints are strict at $p$,
	so they hold in a neighborhood of $p$.
	In the same way, define $T_2\subseteq V$
	from the constraints of $K_2$.
	Since there are finitely many constraints,
	there exists $r>0$ such that, for
	$\|u\|_2<r$ and $\|v\|_2<r$, with $v\in V$,
	\begin{equation}\label{eq:tangent-exact}
		p+u\in K_1\iff u\in T_1,\qquad
		p+v\in K_2\iff v\in T_2.
	\end{equation}
	
	Let $\Gamma_p:=\{g\in\Gamma:U_gp=p\}$.
	For $g\in\Gamma_p$ we have
	$U_g(p+v)-(p+v)=U_gv-v$ and
	$\inner{p}{U_gv-v}=0$. Therefore
	\begin{equation}\label{eq:dirichlet-exact}
		\inner{p+u}{U_g(p+v)-(p+v)}
		=\inner{u}{U_gv-v}
		\qquad(g\in\Gamma_p).
	\end{equation}
	For $g\notin\Gamma_p$, we have
	\[
	\inner{p}{U_gp-p}
	=-\frac12\|U_gp-p\|_2^2<0.
	\]
	Since $\Gamma$ is finite, these constraints
	remain strict near $(p,p)$.
	
	Let $\Psi$ be the feasible set obtained from
	$\mathcal F$ by replacing $K_1,K_2,\Gamma$
	with $T_1,T_2,\Gamma_p$, respectively.
	By~\eqref{eq:tangent-exact} and
	\eqref{eq:dirichlet-exact}, we have
	\[
	(p+u,p+v)\in\mathcal F
	\iff (u,v)\in\Psi
	\]
	whenever $\|u\|_2,\|v\|_2<r$.
	
	We now show that $\inf_\Psi\rho=\nu$.
	For sufficiently large $n$,
	$(u_n-p,v_n-p)\in\Psi$. Since
	$B_1p=B_2p$, we have
	\[
	\rho(u_n-p,v_n-p)=\rho(u_n,v_n).
	\]
	Thus $\Psi\neq\emptyset$ and
	$\inf_\Psi\rho\le\nu$.
	
	Conversely, let $(u,v)\in\Psi$.
	Since $T_1$ and $T_2$ are cones,
	$(su,sv)\in\Psi$ for every $s>0$.
	For sufficiently small $s>0$, the local
	equivalence above gives
	$(p+su,p+sv)\in\mathcal F$.
	Moreover,
	\[
	\rho(p+su,p+sv)=\rho(u,v).
	\]
	Hence $\nu\le\inf_\Psi\rho$.
	Therefore $\inf_\Psi\rho=\nu$, and
	any minimizer on $\Psi$ gives a
	minimizer on $\mathcal F$.
	
	It remains to show that the infimum on
	$\Psi$ is attained.
	Since $\inner{a_i}{p}=0$ for the
	constraints defining $T_1$, we have
	$T_1+\R p=T_1$.
	Similarly, $T_2+\R p=T_2$, since $p\in V$.
	Also, $U_gp=p$ for $g\in\Gamma_p$ and
	$B_1p=B_2p$. Thus $\Psi$ and $\rho$
	are invariant under
	\[
	(u,v)\mapsto(u+tp,v+tp),
	\qquad t\in\R.
	\]
	Taking
	\[
	t=-\frac{\inner{v}{p}}{\|p\|_2^2},
	\]
	we can move the second entry into
	$p^\perp$ without changing feasibility
	or the value of $\rho$.

Let $V':=V\cap p^\perp$, and let
\[
\Psi':=\bigl\{(u,v)\in T_1\times(T_2\cap p^\perp):
u\neq v,\ \inner{u}{U_gv-v}\le0
\ \text{for all }g\in\Gamma_p\bigr\}.
\]
Then $\Psi'\neq\emptyset$ and
$\inf_\Psi\rho=\inf_{\Psi'}\rho$.
Notice that $V'$ is a linear subspace,
$T_1$ and $T_2\cap p^\perp$ are closed
convex polyhedral cones, and $\Gamma_p$
is a finite orthogonal group.
Since $p\in V$ and $p\neq0$, we have
$\dim V'=\dim V-1$.
Thus, by the induction hypothesis,
the infimum on $\Psi'$ is attained.
The same is therefore true for $\Psi$,
and hence for $\mathcal F$.

This completes the proof.
\end{proof}

\begin{cor}\label{cor:main-from-formula}
	Let $f$ be as in Definition~\ref{piecelinear}. Then $a(f)$ is a minimum: there are
	$x_{*},y_{*}\in\R^{d}$ with $[x_{*}]\neq[y_{*}]$ and
	$\|f(x_{*})-f(y_{*})\|_2=a(f)\,d_{G}(x_{*},y_{*})$. In particular, if $\widehat f$ is
	injective on $\R^{d}/G$ then $a(f)>0$. Together with
	Proposition~\ref{prop:upper-lipschitz} this gives a different proof of
	Theorem~\ref{thm:main}, through the optimal lower constant.
\end{cor}

\begin{proof}
	By Proposition~\ref{thm:formula}(i), we have $a(f)=\nu(\alpha_{0},\beta_{0})$ for some pair
	$(\alpha_{0},\beta_{0})\in I\times I$, since $I\times I$ is finite. By
	Proposition~\ref{prop:mu-attained}, the infimum defining $\nu(\alpha_{0},\beta_{0})$ is a
	minimum, taken at some $(x_{*},y_{*})$ with $x_{*}\in C_{\alpha_{0}}$, $y_{*}\in C_{\beta_{0}}$,
	$x_{*}\neq y_{*}$ and $\|x_{*}-y_{*}\|_2=d_{G}(x_{*},y_{*})$. Then $d_{G}(x_{*},y_{*})>0$, so
	$[x_{*}]\neq[y_{*}]$, and the ratio at this pair is $\nu(\alpha_{0},\beta_{0})=a(f)$. Notice that
	$f$ restricts to $A_{\alpha_{0}}$ on $C_{\alpha_{0}}$ and to $A_{\beta_{0}}$ on $C_{\beta_{0}}$.
	Therefore
	$\|f(x_{*})-f(y_{*})\|_2=\|A_{\alpha_{0}}x_{*}-A_{\beta_{0}}y_{*}\|_2=a(f)\,d_{G}(x_{*},y_{*})$.
	If $\widehat f$ is injective, then $f(x_{*})\neq f(y_{*})$, and thus $a(f)>0$. Finally,
	$a(f)\leq b(f)$, since $a(f)$ is an infimum of ratios, each of which is at most $b(f)$.
	Moreover, $b(f)\leq L_f<\infty$ by Proposition~\ref{prop:upper-lipschitz}. Consequently the
	inequalities of Theorem~\ref{thm:main} hold with $a=a(f)$ and $b=b(f)$ for all $x,y$ with
	$[x]\neq[y]$. For $[x]=[y]$ all three terms vanish, because $d_{G}(x,y)=0$ and $f(x)=f(y)$ by
	$G$-invariance.
\end{proof}

Note that Proposition~\ref{prop:mu-attained} is used only in this corollary. In particular,
Proposition~\ref{thm:formula} does not use it.

\section{Numerical Experiments}\label{sec:experiments}

\begin{figure}[!htb]
	\centering
	\includegraphics[width=0.95\linewidth]{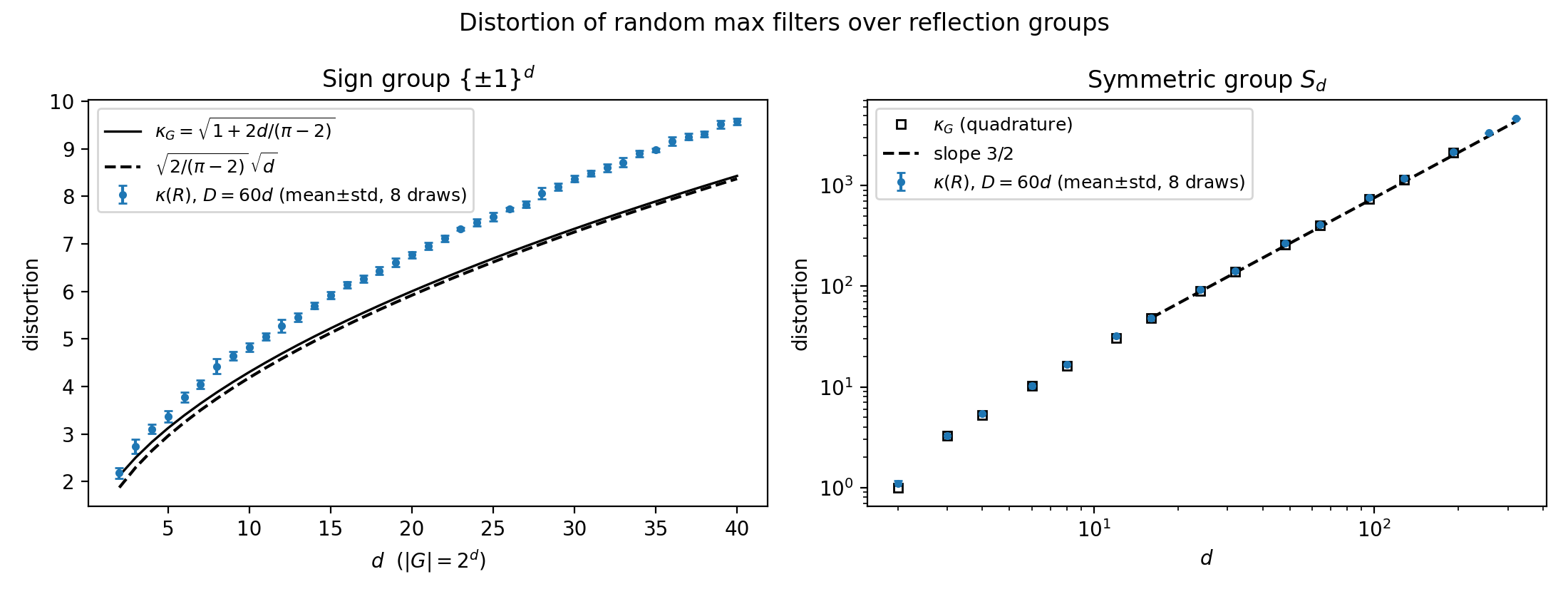}
	\caption{Distortion $\kappa(f_A)=\kappa(R)$ of the random $G$-invariant max filter,
		where $R$ is the matrix of Theorem~\ref{thm:reflection-constants}. The dots are sampled
		at $D=60d$, with mean $\pm$ standard deviation over $8$ template draws. Left: the sign group
		$\{\pm1\}^d$ for $2\le d\le40$. The solid curve is the exact $D\to\infty$ limit
		$\kappa_G=\sqrt{1+2d/(\pi-2)}$ of Theorem~\ref{thm:exact}(a), and the
		dashed curve is its asymptote $\sqrt{2/(\pi-2)}\sqrt d$ as $d\to\infty$. The sampled values lie above
		$\kappa_G$. Right: $S_d$ for $2\le d\le320$, on log-log axes. The squares are the values of $\kappa_G$ computed
		by quadrature for $d\le192$, and the dashed line has slope $3/2$ and passes through the sampled value
		at $d=16$.}
	\label{fig:reflection}
\end{figure}

In this section we test numerically the growth laws of Theorems~\ref{thm:exact} and~\ref{thm:growth}.\footnote{The experiments can be reproduced with the scripts \texttt{verify\_distortion.py} (Figure~\ref{fig:reflection}, sampled values) and \texttt{quadrature\_MG.py} (the quadrature values of $M_G$), which share a helper module. They use \texttt{numpy} and \texttt{scipy}, and \texttt{matplotlib} for the figure.}

Recall that by Theorem~\ref{thm:reflection-constants} the distortion is the condition number of $R$. For small $d$ we checked  $f_A=R\,\Lambda_G$ against a direct evaluation of the max filter, and the two agree. We estimate the distortion in two ways. First, we sample $R$ at a finite $D$. Second, we compute $M_G$ by Gauss Jacobi quadrature on the order statistic moments, and by Theorem~\ref{thm:spectral-limit} this gives the limit $D\to\infty$. We use $220$ quadrature nodes. If we double the number of nodes, the values of $d^2\lambda_{\min}(M_G)$ and $\kappa_G/d^{3/2}$ change by less than $10^{-5}$. Moreover, the identity $\operatorname{tr}M_G=d$ holds to relative accuracy $4\cdot10^{-5}$. The results are shown in Figure~\ref{fig:reflection}.

For the sign group, Theorem~\ref{thm:exact}(a) gives $\kappa_G/\sqrt d\to\sqrt{2/(\pi-2)}\approx1.32$. With $D=60d$ and eight template draws, the measured $\kappa/\sqrt d$ stays between $1.50$ and $1.58$ for $2\le d\le40$. These values are larger than $\kappa_G/\sqrt d$ at the same $d$ by a factor of about $1.03$ to $1.14$. For large $d$ this is consistent with the finite-$D$ inflation factor $(1-\sqrt{d/D})^{-1}\approx1.15$.

For $S_d$, Theorem~\ref{thm:growth}(a) gives the asymptotic upper bound $\sqrt{2/\pi}\approx0.798$ for $\kappa_G/d^{3/2}$. The quadrature values at $d=64$, $96$, $128$ and $192$ are $0.785$, $0.789$, $0.791$ and $0.793$, so they increase toward this bound. At $d=128$ the sampled value at $D=60d$ is about $3\%$ higher. At the same four values of $d$, $d^2\lambda_{\min}(M_G)$ equals $1.5670$, $1.5681$, $1.5688$ and $1.5694$, while $\pi/2\approx1.5708$. This suggests that the floor of Lemma~\ref{lem:orderstat-floor} is asymptotically sharp, as the classical first-order approximation~\cite{stephens1990asymptotic},~\cite[Ch.~4]{david2003order} predicts. However, we do not prove this.

For $B_d$ and $D_d$ at $d=192$, the values of $d^2\lambda_{\min}(M_G)$ are $0.3909$ and $0.3914$. Both are just below $\pi/8\approx0.3927$, which is the limit of $d^2$ times the proved floor. Moreover, the values of $\kappa_G/d^{3/2}$ are $1.5967$ and $1.5957$, and both are below the bound $\sqrt{8/\pi}\,\sqrt{(1+1/d)(1+2/d)}\approx1.6082$ of Theorem~\ref{thm:growth}(b). Finally, for $I_2(k)$ with $2\le k\le32$, the sampled spectrum agrees with the closed form $\{1\pm\sin w/w\}$, $w=\pi/k$, of Theorem~\ref{thm:exact}(b) to within $0.12\%$ relative error.
\FloatBarrier

\section*{Acknowledgments}

I thank Radu Balan and Nadav Dym for reading an earlier version of this manuscript and for their helpful comments and suggestions. I am particularly grateful to Radu Balan for his extensive feedback, which substantially improved the presentation of this work.

\paragraph{Use of AI tools.}
Generative AI tools were used to help identify potentially relevant literature and improve the wording of the manuscript. The author is responsible for verifying the references and for the mathematical content and conclusions.

\printbibliography

\end{document}